\documentclass[12pt]{article}
\usepackage{amsmath}
\usepackage{pictex}
\usepackage{amscd}
\usepackage{amstext}
\usepackage{amssymb}
\usepackage{amsfonts}
\usepackage{latexsym}
\usepackage{epsfig}
\usepackage{verbatim}
\usepackage[dvips]{color}
\title{Complex Interpolation  of  Weighted Besov- and Lizorkin-Triebel Spaces (long version)}

\author{Winfried Sickel\footnote{Friedrich-Schiller-Universit\"at Jena, Mathematisches Institut, Ernst-Abbe-Platz 2, 07743 Jena, Germany, e-mail: {\tt winfried.sickel@uni-jena.de}},
Leszek Skrzypczak\footnote{Adam Mickiewicz University Poznan, Faculty of Mathematics and Computer Science, Ul. Umultowska 87, 61-614 Pozna\'n, Poland, e-mail: {\tt lskrzyp@amu.edu.pl}},
and Jan Vyb\'\i ral\footnote{Technical University Berlin, Department of Mathematics, Secretary office MA 4-1, Street of 17. June 136, 10623 Berlin, Germany, e-mail: {\tt vybiral@math.tu-berlin.de}}
\thanks{Research supported by the DFG Research Center {\sc Matheon} ``Mathematics for key technologies'' in Berlin.}}
\date{\today}

\newtheorem{T}{Theorem}
\newtheorem{Lem}[T]{Lemma}
\newtheorem{Prop}[T]{Proposition}
\newtheorem{Def}[T]{Definition}
\newtheorem{Cor}[T]{Corollary}
\newtheorem{Rem}{Remark}

\newcommand{\bt}{\begin{T}}
\newcommand{\bl}{\begin{Lem}}
\newcommand{\bp}{\begin{Prop}}
\newcommand{\bc}{\begin{Cor}}
\newcommand{\bd}{\begin{Def}}
\newcommand{\br}[2]{\begin{Rem}\label{#1}{\rm #2}}
\newcommand{\er}{ \end{Rem}}
\newcommand{\et}{\end{T}}
\newcommand{\el}{\end{Lem}}
\newcommand{\ep}{\end{Prop}}
\newcommand{\ec}{\end{Cor}}
\newcommand{\ed}{\end{Def}}

\newcommand{\be}{\begin{equation}}
\newcommand{\ee}{\end{equation}}
\newcommand{\beq}{\begin{eqnarray}}
\newcommand{\eeq}{\end{eqnarray}}
\newcommand{\beqq}{\begin{eqnarray*}}
\newcommand{\eeqq}{\end{eqnarray*}}

\def\lsim{\raisebox{-1ex}{$~\stackrel{\textstyle <}{\sim}~$}}

\def\gsim{\raisebox{-1ex}{$~\stackrel{\textstyle >}{\sim}~$}}

\def\tr{\mathop{\rm tr \,}\nolimits}
\def\ext{\mathop{\rm ext \,}\nolimits}

\newcommand{\bpr}{{\bf Proof.}\ }
\newcommand{\epr}{\hspace*{\fill}\rule{3mm}{3mm}\\}

\newcommand{\bspq}{B^s_{p,q}}

\newcommand{\fspq}{F^s_{p,q}}

\newcommand{\Rd}{{\mathbb R}^{d}}
\newcommand{\R}{{\mathbb R}}
\newcommand{\N}{{\mathbb N}}
\newcommand{\Z}{\mathbb Z}
\newcommand{\Zd}{\Z^{d}}
\newcommand{\C}{{\mathbb C}}
\newcommand{\cs}{{\mathcal{S}}}
\newcommand{\ca}{{\mathcal{A}}}
\newcommand{\cx}{{\mathcal{X}}}
\newcommand{\cl}{{\mathcal{L}}}

\newcommand{\aspq}{A^s_{p,q}}

\newcommand{\cf}{{\mathcal{F}}\,}
\newcommand{\cfi}{{\mathcal{F}}^{-1}\,}

\newlength{\fixboxwidth}
\begin{document}

\maketitle

\begin{abstract}
We study complex interpolation of weighted Besov and Lizorkin-Triebel spaces.
The used weights $w_0,w_1$ are local Muckenhoupt weights in the sense of Rychkov.
As a first step we calculate the  Calder\'on products of associated sequence spaces.
Finally, as a corollary  of these investigations, we obtain results on  complex interpolation of
radial subspaces of Besov and Lizorkin-Triebel spaces on $\R^d$. 
\end{abstract}



\newpage

\section{Introduction}


Nowadays interpolation theory has been established as an
important tool in various branches of mathematics, in particular in analysis of PDE's.
Within the known interpolation methods the complex interpolation method of Calder{\'o}n, denoted by $[\, \cdot\, , \, \cdot\,]_\Theta$,
is of particular importance and probably the most often used one.\\
Let $L_{p}(\Rd,w)$ denote the weighted Lebesgue space with weight $w$.
Here in this paper we study
generalizations of the following formula
\beq\label{wichtig5}
\Big[L_{p_0} (\Rd,w_0), L_{p_1}(\Rd,w_1)\Big]_\Theta & = &  L_{p}(\Rd,w)\, , \qquad 1\le p_0,p_1 < \infty\, ,
\\
\label{wichtig6}
\qquad \frac 1p & := &  \frac{1-\Theta}{p_0} + \frac{\Theta}{p_1}\, , \qquad w  :=  w^{(1-\Theta)\, p/p_0}_0\, w^{\Theta\,  p/p_1}_1\, ,
\eeq
see, e.g., \cite[Theorem.~5.5.3]{BL} or \cite[Theorem~1.18.5]{TrI}.
We shall replace the weighted Lebesgue spaces $L_{p}(\Rd,w)$
by weighted Besov and Lizorkin-Triebel spaces.
There are already some contributions dealing with this problem. Let us mention here
Bownik \cite{Bow2} and Wojciechowska \cite{AW2}. But both authors only deal with the case
$w=w_0=w_1$.
Bownik considers weights related to doubling measures and
Wojciechowska is dealing with local Muckenhoupt weights (as we shall do in most of the cases).
Whereas in (\ref{wichtig5}) it will be enough that $w_0$ and $w_1$ are positive, in the generalizations,
we have in mind, it is not clear what is the correct class of weights.
It seems that necessary conditions concerning the weights are not known in this context.
\\
To calculate
\be\label{wichtig7}
\Big[F^{s_0}_{p_0,q_0}(\R^d,w_{0}), \, F^{s_1}_{p_1,q_1}(\R^d,w_{1})\Big]_\Theta  \qquad \mbox{and}\qquad
\Big[B^{s_0}_{p_0,q_0}(\R^d,w_{0}), \, B^{s_1}_{p_1,q_1}(\R^d,w_{1})\Big]_\Theta
\ee
we shall apply a method which has been used by Bownik \cite{Bow2} and
Wojciechowska \cite{AW2} as well.
First we shall calculate the Calder{\'o}n products of associated sequence spaces.
Afterwards we shall use the known coincidence of
Calder{\'o}n products and the complex method of interpolation (under certain extra conditions)
to shift the results to the complex interpolation of these sequence spaces.
Finally, these results are lifted by wavelet isomorphisms to the level of function spaces.
This method has been developed by Calder\'on \cite{Cal},
Frazier, Jawerth \cite{FJ2},  Mendez, Mitrea \cite{MM} and  Kalton, Maybororda, Mitrea \cite{KMM}.
The latter two references are connected with the extension of the complex method to certain quasi-Banach spaces.
Also in our paper we shall work with quasi-Banach spaces.
In fact,  we will allow the maximal range of the parameters in (\ref{wichtig7}) with the exception of
$F^{s}_{\infty,q}(\R^d,w)$. Of course, it would be interesting to incorporate these spaces as well but this requires additional effort.
\\
The paper is organized as follows.
In Section \ref{real} we recall the basic notions for the complex method
and also describe the state of the art in the unweighted case.
The next section is devoted to the calculation of the
Calder\'on products of some weighted sequence spaces. For  $w_0, w_1 \in {\mathcal A}_\infty^{\ell oc}$
we define
\[
w(x) := w_{0} (x)^{\frac{(1-\Theta)p}{p_0}}\, w_1(x)^{\frac{\Theta p}{p_1}}\, ,
\qquad x \in \Rd\, .
\]
Then we will  establish the formulas
\[
f^{s_0}_{p_0,q_0}(\R^d,w_0)^{1-\Theta}\,  f^{s_1}_{p_1,q_1}(\R^d,w_1)^\Theta  =
f^{s}_{p,q}(\R^d,w)
\]
as well as
\[
b^{s_0}_{p_0,q_0}(\R^d,w_0)^{1-\Theta}\,  b^{s_1}_{p_1,q_1}(\R^d,w_1)^\Theta  =
b^{s}_{p,q}(\R^d,w)
\]
under rather general conditions on the parameters.
This will be the most complicated part.
In Section \ref{Complex1} we deal with
complex interpolation of weighted Besov and
Lizorkin-Triebel spaces.
On the one side we simply shift here the results, obtained for Calder{\'o}n products, to complex interpolation formulas,
on the other side we apply the result of Shestakov \cite{Sh} to calculate
(\ref{wichtig7}) also in some of those situations where both spaces are not separable.
Finally, in Section \ref{complex} we apply the results obtained before to derive some
complex interpolation formulas for radial subspaces of Besov and
Lizorkin-Triebel spaces.
This was actually the original motivation for our work.
\\
Definitions and some properties of the classes of weights and classes of function spaces under consideration here
are collected in the Appendix at the end of this paper.


\subsection*{Notation}


As usual, $\N$ denotes the natural numbers, $\N_0:= \N \cup \{0\}$,  $\Z$ denotes the integers and $\R$ the real numbers.
For the  complex numbers we use the symbol $\C$, for the Euclidean $d$-space  we use $\Rd$
and $\Zd$ to denote the collection of all elements in $\Rd$ having integer components.
At very few places we shall need the Fourier transform $\cf$ as well as its inverse transformation
$\cfi$, always defined on the Schwartz space  $\cs' (\Rd)$ of tempered distributions.
\\
If $X$ and $Y$ are two quasi-Banach spaces, then the symbol  $X \hookrightarrow Y$
indicates that the embedding is continuous.
As usual, the symbol  $c $ denotes positive constants
which depend only on the fixed parameters $s,p,q$ and probably on auxiliary functions,
unless otherwise stated; its value  may vary from line to line.
Sometimes we  will use the symbols ``$ \lsim $''
and ``$ \gsim $'' instead of ``$ \le $'' and ``$ \ge $'', respectively. The meaning of $A \lsim B$ is given by: there exists a constant $c>0$ such that
 $A \le c \,B$. Similarly $\gsim$ is defined. The symbol
$A \asymp B$ will be used as an abbreviation of
$A \lsim B \lsim A$.\\
Inhomogeneous weighted Besov and Lizorkin-Triebel spaces are denoted by
$B^{s}_{p,q}(\R^d,w)$ and $F^{s}_{p,q}(\R^d,w)$, respectively.
If the weight is identically $1$, we shall drop $w$ in notation.
In case  there is no reason to distinguish between these two scales we will use the
notation   $A^{s}_{p,q}(\R^d,w)$.
Definitions, properties  as well as some references are given in the Appendix.


\section{Complex interpolation of  Besov and
Lizorkin-Triebel spaces - the state of the art}
\label{real}


For convenience of the reader we recall some notions from interpolation theory
as well as some results in the framework of
Besov and Lizorkin-Triebel spaces. \\
For the basics of interpolation theory we refer to the monographs
\cite{BL,TrI,KPS,BS}.


\subsection{The complex method of interpolation}


The complex method in case of interpolation couples of Banach spaces is a well-studied subject,
see the quoted monographs above.
Here we are interested in the complex method in case of interpolation couples of certain quasi-Banach spaces.
For that reason we give some details.
We follow \cite{KMM}, see also \cite{MM} and \cite{KM}.

\bd\label{Aconvex}
A quasi-Banach space $(X,\|\cdot|X\|)$ is called {\it analytically convex} if there
is a constant $c$ such that for every polynomial $P:\C\to X$ we have
$$\|P(0)|X\|\le c\max_{|z|=1}\|P(z)|X\|.$$
\ed

In the framework of  analytically convex quasi-Banach spaces
the Maximum Modulus Principle holds.
Let
\[
S_0:= \{z \in \C: \: 0<\Re e\, z<1\}\qquad \mbox{and} \qquad S:=\{z\in\C:\: 0\le \Re e\, z\le 1\}\, .
\]

\begin{Prop} For a quasi-Banach space $(X,\|\cdot|X\|)$ the following
conditions are equivalent:

{\rm(i)} X is analytically convex.

{\rm(ii)} There exists a constant  $c$ such that
$$
\max\{\|f(z)|X\|:z\in S_0\}\: \le\:  c\, \max\{\|f(z)|X\|:z\in S\setminus S_0\}
$$
for any function $f:S\to X$, analytic on $S_0$ and continuous and bounded on
$S$.
\end{Prop}

We refer to \cite[Theorem~7.4]{KMM}.
Based on this property the following definition makes sense.

\begin{Def}
Let $(X_0,X_1)$ be an interpolation couple of quasi-Banach spaces,
i.e., $X_0$ and $X_1$ are continuously embedded into a larger topological vector space $Y$. 
In addition, let $X_0 + X_1$ be analytically convex.
Let ${\mathcal A}$ be the set of all bounded and analytic functions $f: S_0 \to X_0 + X_1$, which
extend continuously to the closure $S$ of the strip s.t. the traces $t \mapsto f(j+it)$ are bounded continuous
functions into $X_j$, $j=0,1$. We endow $\ca$ with the quasi-norm
\[
\|\, f \, |\ca \| := \max \Big\{\sup_{t\in \R}\, \|\, f(it)\, |X_0\|, \,  \sup_{t\in \R}\, \|\, f(1+it)\, |X_1\|\Big\}\, .
\]
Let $0 < \Theta < 1$.
Further, we define $[X_0,X_1]_\Theta$ to be the set of all $x \in \ca (\Theta):= \{f(\Theta): \, f \in \ca \}$ and
\[
\|\, x \, |[X_0,X_1]_\Theta\| := \inf \Big\{\|\, f\, |\ca\|: \: f(\Theta) = x\Big\}\, .
\]
\end{Def}

\begin{Rem}
 \rm
Any Banach space is analytically convex.
Hence, if $(X_0,X_1)$ is an interpolation couple of Banach spaces,
this reduces to the standard definition of $[X_0,X_1]_\Theta$.
\end{Rem}

\begin{Lem}\label{convex}
Let $0 < q \le \infty$, $s \in \R$ and $w \in A_\infty^{\ell oc}$.
\\
(i) Let $0 < p< \infty$. Then
$F^s_{p,q} (\Rd,w)$ is analytically convex.
 \\
(ii) Let $0 < p \le  \infty$. Then
$B^s_{p,q} (\Rd,w)$ is analytically convex.
\end{Lem}

\noindent
\bpr
If $X$ is an analytically convex quasi-Banach space and $Y$ is a closed subspace of $X$ then
 $Y$ is analytically convex, see Proposition~7.5 in \cite{KMM}.
By means of Proposition \ref{sawano} it will be enough to prove analytic convexity for
the sequence spaces $f^s_{p,q} (\Rd,w)$  and $b^s_{p,q} (\Rd,w)$.
In contrast to the function spaces the sequence spaces are quasi-Banach lattices.
We need a further notion.
A quasi-Banach lattice of functions $(X,\|\cdot|X\|)$ is called
{\it lattice r-convex} if
$$
\Bigl\|\bigl(\sum_{j=1}^m|f_j|^r\bigr)^{1/r}|X\Bigr\|\le
\Bigl(\sum_{j=1}^m\|f_j|X\|^r\Bigr)^{1/r}
$$
for any finite family $\{f_j\}_{1\le j\le m}$ of functions from $X$.
\\
There are simple criteria for a quasi-Banach lattice of functions
to be analytically convex, see \cite[Theorem~7.8]{KMM}:
$X$ is analytically convex if, and only if, $X$ is lattice $r$-convex for some $r>0$.
\\
It remains to show that  the sequence spaces $a^s_{p,q} (\Rd,w)$, $a \in \{b,f\}$, are lattice $r$-convex.
This holds with $r \le \min (p,q,1)$ by standard arguments (use the generalized Minkowski inequality).
\epr

\begin{Rem}
 \rm
The unweighted case was considered in Mendez and Mitrea \cite{MM}, see also Kalton, Mayboroda and Mitrea
\cite[Proposition~7.7]{KMM}.
For weights related to doubling measures the statement has been settled by Bownik \cite{Bow2}.
\end{Rem}


\subsection{The state of the art}


For convenience of the reader we recall what is known in the unweighted
situation. Comments to the weighted case will be given within the text.

\begin{Prop}\label{interpol2}
Let  $0 <  p_0 ,  p_1 \le \infty$, $0 <  q_0 ,\,  q_1 \le \infty$, $s_0, s_1 \in \R$, and $0 < \Theta < 1$.
Define
\be\label{def}
s:= (1-\Theta)\, s_0 + \Theta \, s_1\, , \qquad
\frac 1p := \frac{1-\Theta}{p_0} + \frac{\Theta}{p_1}
\qquad \mbox{and}\qquad \frac 1q := \frac{1-\Theta}{q_0} +
\frac{\Theta}{q_1} \, .
\ee
{\rm (i)} Let $\min (p_0+q_0,p_1+q_1) <\infty$.
Then we have
\be\label{eq-30}
B^s_{p,q} (\Rd)  =   \Big[B^{s_0}_{p_0,q_0} (\Rd),
B^{s_1}_{p_1,q_1} (\Rd)\Big]_{\Theta} \, .
\ee
{\rm (ii)} Let $\max(p_0,p_1) <\infty$ and $\min (q_0,q_1) <\infty$.
Then we have
\be\label{eq-29}
F^s_{p,q} (\Rd) =   \Big[F^{s_0}_{p_0,q_0} (\Rd),
F^{s_1}_{p_1,q_1} (\Rd)\Big]_{\Theta} \, .
\ee
\end{Prop}

\begin{Rem}\label{Grenzfaelle}
 \rm
(i) Proposition \ref{interpol2} in this generality can be found in Frazier and Jawerth \cite{FJ2}
(the F-case) and in Kalton, Mayboroda, Mitrea \cite[Theorem~9.1]{KMM}(F- and  B-case).
With the extra condition $s_0 \neq s_1$ one can find
(\ref{eq-30}) also in Mendez, Mitrea \cite{MM}.
However, Proposition \ref{interpol2}  has many forerunners in case $\min (p_0,p_1,q_0, q_1)\ge 1$, e.g.,
Calder{\'o}n, J.L.~Lions, Magenes, Taibleson, Grisvard, Schechter, Peetre and Triebel.
We refer to \cite{BL,Pe,TrI} and the references given there.
\\
(ii)
Let us mention, that formula (\ref{eq-29}) remains true in case if either $\max (p_0,q_0) <  \infty$
or $\max (p_1,q_1) <  \infty$, see \cite{FJ2} and \cite{KMM}.
However, in our paper we shall not deal with
the spaces $F^s_{\infty,q} (\Rd)$ and its weighted counterparts.\\
(iii)
The counterpart of (\ref{eq-29}) for anisotropic Lizorkin-Triebel spaces
(more exactly, the generalization to)
has been proved by Bownik \cite{Bow2}.
\\
(iv) It is of certain interest to notice
that in some cases complex interpolation of pairs of Besov spaces
does not result in a Besov space.
More exactly, if $1 < p< \infty$, $s_0,s_1 \in \R$, $s_0 \neq s_1$,
$s:= (1-\Theta)\, s_0 + \Theta \, s_1$,
then
\be\label{eq-30not}
\mathring{B}^s_{p,\infty} (\Rd)  =   \Big[B^{s_0}_{p,\infty} (\Rd),
B^{s_1}_{p,\infty} (\Rd)\Big]_{\Theta} \, ,
\ee
where $\mathring{B}^s_{p,\infty} (\Rd)$ denotes the closure of the set of test functions in
${B}^s_{p,\infty} (\Rd)$, a space strictly smaller than ${B}^s_{p,\infty} (\Rd)$.
We refer to \cite[Theorem~2.4.1]{TrI}. Later on, see Subsection \ref{gap}, we shall supplement this formula.
\\
(v) In \cite[Theorem~6.4.5]{BL} the following formula is claimed to be true:
\be\label{eq-30notb}
{B}^s_{p,q} (\Rd)  =   \Big[B^{s_0}_{p_0,q_0} (\Rd),
B^{s_1}_{p_1,q_1} (\Rd)\Big]_{\Theta} \, ,
\ee
where $1 \le p_0,p_1,q_0,q_1\le \infty$, $s_0 \neq s_1$ and $s,p,q$ as in (\ref{def}).
Of course, this is in contradiction with part (iv) in case $q=q_0=q_1=\infty$.
We do not believe in this formula in cases, when  both spaces $B^{s_0}_{p_0,q_0} (\Rd)$ and
$B^{s_1}_{p_1,q_1} (\Rd)$ are not separable, see Subsection \ref{gap} for more information.
\\
(vi) There is a number of further methods of interpolation where
the outcome is known in case of pairs of either Besov or
Lizorkin-Triebel spaces. Most prominent is the real method of interpolation.
For corresponding results we refer to \cite[Theorem~6.4.5]{BL} and \cite[2.4]{Tr83}.
Triebel  \cite[2.4]{Tr83} also had invented a certain modification of the complex method
and has been able to prove the counterparts of (\ref{eq-30}), (\ref{eq-29}) for this modified complex method.
However, it is not known whether this modified method  has the interpolation property.
Frazier, Jawerth \cite{FJ2}  and Bownik \cite{Bow2} also studied the $\pm$-method of Gustaffson and Peetre, denoted by $\langle X_0, X_1,\theta\rangle$,
and the method $\langle X_0, X_1\rangle_\theta$, due to Gagliardo.
Then Proposition \ref{interpol2}(ii) remains true also for these methods.
\end{Rem}


\section{Calder{\'o}n products  of sequence spaces associated to weighted Besov and
Lizorkin-Triebel spaces}
\label{Calderon}


After having introduced the necessary definitions in Subsection \ref{Calderona}
we shall first deal with the Calder{\'o}n products of the sequence spaces $f^s_{p,q}(\Rd,w)$
(originated from the ingenious proof of Frazier and Jawerth in the unweighted situation).
In the third subsection we shall investigate Calder{\'o}n products of the sequence spaces $b^s_{p,q}(\Rd,w)$
by employing a totally different method.


\subsection{Definition and basic properties of the Calder{\'o}n product}
\label{Calderona}


Let $({\mathfrak X}, {\mathcal S}, \mu)$ be a $\sigma-$finite measure space and
let ${\mathfrak M}$ be the class of all complex--valued, $\mu-$measurable functions on
${\mathfrak X}$. Then a quasi-Banach space $X\subset {\mathfrak M}$
is called a {\it quasi-Banach lattice of functions}
 if for every $f\in X$ and $g\in{\mathfrak M}$ with $|g(x)|\le |f(x)|$ for $\mu-$a.e. $x\in{\mathfrak X}$
one has $g\in X$ and $\|g|X\|\le \|f|X\|.$

\begin{Def}
Let $({\mathfrak X}, {\mathcal S}, \mu)$ be a $\sigma-$finite measure space and
let ${\mathfrak M}$ be the class of all complex--valued, $\mu-$measurable functions on
${\mathfrak X}$. Let
$X_j \subset {\mathfrak M}$, $j=0,1$,  be quasi-Banach lattices of functions.
Let $0<\Theta<1$. Then the {\it Calder\'on product $X_0^{1-\Theta}X_1^\Theta$}
of $X_0$ and $X_1$ is the collection of all functions $f \in {\mathfrak M}$ s.t.
the  quasi-norm
\beqq
\|f|X_0^{1-\Theta}X_1^\Theta\| & := & \inf\Bigl\{\|f_0|X_0\|^{1-\Theta}\|f_1|X_1\|^\Theta:\:
|f|\le |f_0|^{1-\Theta}|f_1|^\Theta \quad \mu - \mbox{a.e.},
\\
&& \hspace{7cm} f_j\in X_j, \, j=0,1\Bigr\}
\eeqq
is finite.
\end{Def}

\begin{Rem}
 \rm
(i) Calder{\'on} products have been introduced by Calder{\'o}n \cite[13.5]{Cal}
(in a little bit different form which coincides with the above one).
The usefulness of this method and its limitations  have been perfectly described by Frazier and Jawerth \cite{FJ2}
which we quote now:
{\it Although restricted to the case of a lattice, the Calder{\'o}n product has the advantage of being
defined in the quasi-Banach case, and,
frequently, of being easy to compute. It has the disadvantage that the interpolation property
(i.e., the property that a linear transformation $T$ bounded on $X_0$ and $X_1$ should be bounded on the spaces in between)
is not clear in general.}
\\
(ii) A further remark to the literature. Calder{\'o}n products
are not investigated in the most often quoted books on interpolation theory:
Bergh and L\"ofstr\"om \cite{BL} (except a short remark on page 129), Triebel \cite{TrI}
and Bennett and Sharpley \cite{BS}.
However, in the monographs of Kre{\u{\i}n}, Petunin and Sem\"enov \cite[pp.~242-246]{KPS},
Brudnyi and Kruglyak \cite[4.3]{BK} and in the lecture note of Maligranda
\cite{Ma} a few informations about
Calder{\'on} products can be found, sometimes in the more general framework of
Calder{\'o}n-Lozanovskii constructions. All these references are concerned with Banach spaces.
Since we need this concept in quasi-Banach spaces as well, we refer in addition to Nilsson
\cite{Ni}, Frazier and Jawerth \cite{FJ2},  Kalton and Mitrea \cite{KM}, Mendez and Mitrea \cite{MM}, Kalton, Mayboroda and Mitrea
\cite{KMM} and Yang, Yuan and Zhuo \cite{DWC}.
\end{Rem}

We collect a few useful properties for later use, see  \cite{DWC}.

\begin{Lem}\label{Leb}
Let $({\mathfrak X}, {\mathcal S}, \mu)$ be a $\sigma-$finite measure space and
let ${\mathfrak M}$ be the class of all complex--valued, $\mu-$measurable functions on
${\mathfrak X}$. Let
$X_j \subset {\mathfrak M}$, $j=0,1$,  be quasi-Banach lattices of functions.
Let $0<\Theta<1$.
\\
{\rm (i)}
Then the {\it Calder\'on product $X_0^{1-\Theta}X_1^\Theta$} is a quasi-Banach space.
\\
{\rm (ii)}
Define $\widetilde{X_0^{1-\Theta}X_1^\Theta}$ as the collection of all $f$ s.t.
there exist a positive real number $\lambda$ and elements $g \in X_0$ and $h \in X_1$
satisfying
\[
|f| \le \lambda |g|^{1-\Theta}\, |h|^\Theta\, , \qquad \|\, g\, |X_0\|\le 1\, , \qquad \|\, h\, |X_1\|\le 1\, .
\]
We put
\[
\|\, f \, |\widetilde{X_0^{1-\Theta}X_1^\Theta}\| := \inf \Big\{ \lambda >0: \quad
|f|\le \lambda |g|^{1-\Theta}\, |h|^\Theta\, , \qquad \|\, g\, |X_0\|\le 1\, , \qquad \|\, h\, |X_1\|\le 1\Big\}\, .
\]
Then
$\widetilde{X_0^{1-\Theta}\, X_1^\Theta} = X_0^{1-\Theta}\, X_1^\Theta$
follows with equality of quasi-norms.
\end{Lem}

Here is one well-known example of a Calder{\'o}n product which can be easily calculated.
Let  $w:~ \Rd \to [0,\infty)$ be measurable and   positive a.e..
The weighted Lebesgue space $L_p (\Rd,w)$
is the collection of all measurable functions $f: \, \Rd \to \C$ such that
\[
\| \, f \, |L_p (\Rd,w)\| := \Big(\int_{\Rd} |f(x)|^p \, w(x)\, dx \Big)^{1/p}<\infty \, .
\]
In case $p=\infty$ we shall use the convention  $L_\infty (\Rd,w) := L_\infty (\Rd)$,  i.e., we always take  $w \equiv 1$.

\begin{Lem}\label{beispiel}
 Let $0 < \Theta < 1$,  $0 < p_0 , p_1 \le \infty $ and let  $w_j:~ \Rd \to [0,\infty)$, $j=0,1$, be measurable and  positive a.e..
We define
\be\label{eq:inter:2b}
\frac 1 p := \frac{1-\Theta}{p_0} + \frac{\Theta}{p_1} \qquad \mbox{and}\qquad
w := w_{0}^{\frac{(1-\Theta)p}{p_0}} \, w_1^{\frac{\Theta p}{p_1}}.
\ee
Then
\[
L_{p_0} (\Rd,w_0)^{1-\Theta} \, L_{p_1} (\Rd,w_1)^{\Theta} = L_{p} (\Rd,w)
\]
with coincidence of the quasi-norms.
\end{Lem}

\noindent
\bpr
{\em Step 1.} Let $\max (p_0,p_1)<\infty$.\\
{\em Substep 1.1} We prove $L_{p_0} (\Rd,w_0)^{1-\Theta} \, L_{p_1} (\Rd,w_1)^{\Theta} \subset L_{p} (\Rd,w)$.
Let $f \in L_{p_0} (\Rd,w_0)^{1-\Theta} \, L_{p_1} (\Rd,w_1)^{\Theta}$. Then there exist
$f_0,f_1$ s.t.
\[
|f(x)| \le |f_0(x)|^{1-\Theta}\, |f_1(x)|^\Theta \qquad \mbox{a.e. in} \: \Rd
\]
and $f_j \in L_{p_j} (\Rd,w_j)$, $j=0,1$.
We employ this inequality together with H\"older's inequality and obtain
\beqq
\Big(\int_{\Rd} |f(x)|^{p} && \hspace{-1.0cm} w(x)\, dx \Big)^{1/p}  \le
\Big(\int_{\Rd} |f_0(x)|^{(1-\Theta)p} \, w_{0} (x)^{\frac{(1-\Theta)p}{p_0}}\, |f_1 (x)|^{\Theta p}
\, w_1(x)^{\frac{\Theta p}{p_1}}\, dx\Big)^{1/p}
\\
& \le &
\Big(\int_{\Rd} |f_0(x)|^{p_0} \, w_{0} (x)\, dx \Big)^{(1-\Theta)/p_0}
\, \Big(\int_{\Rd}  |f_1 (x)|^{p_1} \, w_1(x) \, dx \Big)^{\Theta  /p_1}\, .
\eeqq
Hence, $ f\in L_p (\Rd,w)$ and
$\|\, f\, |L_p (\Rd,w)\| \le \|\, f \, |L_{p_0} (\Rd,w_0)^{1-\Theta} \, L_{p_1} (\Rd,w_1)^{\Theta}\|$.
\\
{\em Substep 1.2.} We prove $L_{p} (\Rd,w) \subset L_{p_0} (\Rd,w_0)^{1-\Theta} \, L_{p_1} (\Rd,w_1)^{\Theta}$.
For given $f \in L_{p} (\Rd,w)$ we define
\[
f_0 (x):= |f(x)|^{p/p_0} \, \Big( \frac{w(x)}{w_0 (x)}\Big)^{1/p_0} \qquad \mbox{and}\qquad
f_1 (x):= |f(x)|^{p/p_1} \, \Big( \frac{w(x)}{w_1 (x)}\Big)^{1/p_1}\, .
\]
Then $f_j \in L_{p_j} (\Rd,w_j)$, $j=0,1$, 
which implies that
$f\in L_{p_0} (\Rd,w_0)^{1-\Theta} \, L_{p_1} (\Rd,w_1)^{\Theta}$ and
\beqq
\| \, f\, | L_{p_0} (\Rd,w_0)^{1-\Theta} \, L_{p_1} (\Rd,w_1)^{\Theta}\| & \le & \| \, f\, | L_{p} (\Rd,w)\|^{(1-\Theta)p/p_0}
\, \| \, f\, | L_{p} (\Rd,w)\|^{\Theta p/p_1}
\\
& = & \| \, f\, | L_{p} (\Rd,w)\|\, .
\eeqq
This proves the claim.
\\
{\em Step 2.} Let $\min (p_0,p_1)< \max (p_0,p_1)=\infty$.
We shall concentrate on the case $0 < p_0 < p_1 = \infty$. Then, by our convention, $w_1 := 1$.
The modifications, needed in Substep 1.1, are obvious.
The function $f_1$, used in Substep 1.2, is now given by $f_1 =1$. With this choice the needed arguments are the same.
\\
{\em Step 3.} Let $p_0 = p_1 = \infty$.
The proof of $L_{\infty} (\Rd)^{1-\Theta} \, L_{\infty} (\Rd)^{\Theta} = L_{\infty} (\Rd)$
is obvious.
\epr

\begin{Rem}
\rm
In case $1\le p_0,p_1\le \infty$ this result can be found in \cite[Exercise~4.3.8]{BK}.
For the unweighted case we also refer to \cite[formula 1.6.1 on page 2.4.6]{KPS} and \cite[Exercise~3 on page 179]{Ma}.
\end{Rem}

Weighted $L_p$-spaces are lattice $r$-convex with $r \le \min (1,p)$, hence analytically convex, see
the proof of Lemma \ref{convex} for an explanation of this notion and \cite[Theorem~7.8]{KMM} for a proof.
Hence, complex interpolation of pairs of weighted $L_p$-spaces makes sense.
There are nice  connections between complex interpolation spaces and
the corresponding Calder\'on product, see the original paper of Calder{\'o}n \cite{Cal} or  Theorem 7.9 in \cite{KMM}.

\begin{Prop}\label{Thm4}
Let $({\mathfrak X}, {\mathcal S},\mu)$ be a complete separable metric space, let $\mu$ be a $\sigma-$finite
Borel measure on ${\mathfrak X}$, and let $X_0,X_1$ be a pair of
quasi-Banach lattices of functions on $({\mathfrak X},\mu)$.
Then, if both $X_0$ and $X_1$ are analytically convex and separable, it follows that
$X_0+X_1$ is analytically convex and
\be\label{sl}
[X_0,X_1]_\Theta = X_0^{1-\Theta}\, X_1^\Theta \, ,\qquad  0<\Theta<1.
\ee
\end{Prop}

Lemma \ref{beispiel} and Proposition \ref{Thm4} immediately imply the following extension of (\ref{wichtig5}).

\begin{Cor}
 Let $0 < \Theta < 1$,  $0 < p_0 , p_1 < \infty $ and let  $w_j:~ \Rd \to [0,\infty)$, $j=0,1$, be measurable and  positive a.e..
Let $p$ and $w$ be defined as in (\ref{eq:inter:2b}). Then
\[
\Big[L_{p_0} (\Rd,w_0), \, L_{p_1} (\Rd,w_1)\Big]_{\Theta} = L_{p} (\Rd,w)
\]
in the sense of equivalence  of  quasi-norms.
\end{Cor}

\begin{Rem}
 \rm
Also in Gustavsson \cite{Gu} and Nilsson \cite{Ni}
 interpolation of $L_{p_0} (\Rd,w_0)$ and $L_{p_1} (\Rd,w_1)$
is discussed  for the full range of $p_0$ and $p_1$.
They considered $\langle L_{p_0} (\Rd,w_0), \, L_{p_1} (\Rd,w_1) \rangle$, where
$\langle \, \cdot \, , \, \cdot\,  \rangle$
denotes an interpolation method introduced by Gagliardo
and coincides with the Calder{\'o}n product under certain conditions, see \cite{Ni}.
\end{Rem}


\subsection{Calder{\'o}n products of  $f^s_{p,q}(\Rd,w)$ spaces}
\label{Calderonb}


In case of weighted  Besov or Lizorkin-Triebel spaces there exist wavelet isomorphisms
which relate these spaces to weighted sequence spaces, see the Appendix for more details.
We first study Calder{\'o}n products of these sequence spaces.\\
Here we are going to use the following abbreviations.
By
\[
Q_{j,k} := \{x \in \Rd: \quad 2^{-j}k_\ell \le x_\ell< 2^{-j}(k_\ell+1)\, , \quad \ell=1, \ldots \, , d\}\, , \quad
j \in \N_0\, , \: k \in \Zd\, ,
\]
we denote the dyadic cubes in $\Rd$ (with volume $\le 1$).
The symbol $\cx_{j,k}$ is used for the characteristic function of the cube $Q_{j,k}$.

\begin{Def}
Let $0<q\le \infty$, $s\in \R$ and let  $w:~ \Rd \to [0,\infty)$ be a  nonnegative measurable function.
In case $0 < p< \infty$ we define
\begin{multline}
f^s_{p,q}(\Rd,w)  :=   \bigg\{ \{\lambda_{j,k}\}_{j,k} : \quad\lambda_{j,k} \in \C\, ,
\\
\| \, (\lambda_{j,k}) \, |f^s_{p,q}(\R^d, w)\| :=
\bigg\|\bigg( \sum_{j=0}^\infty  \sum_{k\in\Z^d}
 2^{sjq}\, |\lambda_{j,k}|^q\, \cx_{j,k}(\cdot)\bigg)^{1/q}\bigg| L_p(\R^d,w)
\bigg\|  < \infty \Bigg\}\, .
\end{multline}
\end{Def}

\noindent
Obviously the spaces $f^s_{p,q}(\Rd,w) $ are quasi-Banach lattices.
For us only those weights $w$ will be of interest which are locally integrable and  satisfy
\be\label{gewicht}
0 < w(Q_{j,k}) := \int_{Q_{j,k}} w(x)\, dx <\infty \qquad \mbox{for all}\quad j\in \N_0, \: k \in \Zd\, .
\ee

\begin{Rem}
\rm
(i) In case $w(x) = 1$ for all $x \in \Rd$ we are back in the unweighted situation.
The associated sequence spaces are denoted simply by
$f^s_{p,q}(\R^d)$.
\\
(ii)
Let $w$ satisfy (\ref{gewicht}).
Let $\mathring{f}^s_{p,q} (\Rd,w)$ denote the closure of the finite sequences in ${f}^s_{p,q} (\Rd,w)$.
Then
\[
\mathring{f}^s_{p,q} (\Rd,w) = {f}^s_{p,q} (\Rd,w) \qquad \Longleftrightarrow \qquad q < \infty.
\]
Especially, if $q = \infty$, then $\mathring{f}^s_{p,\infty} (\Rd,w)$ is a proper subspace of
${f}^s_{p,\infty} (\Rd,w)$.
\\
(iii) Let $w$ satisfy (\ref{gewicht}). It is easily checked that $f^s_{p,q}(\Rd,w)$ is separable if, and only if,
$q < \infty$.\\
(iv) Frazier and Jawerth have introduced also the spaces ${f}^s_{\infty,q} (\Rd)$.
Here in this paper we shall not deal with these classes.
\end{Rem}

Now we turn to the investigation of the Calder{\'o}n products of these sequence spaces.
As mentioned above we are interested in the most general situation
(except the use of $f^s_{\infty,q} (\Rd,w)$). The class of weights, we are dealing with, is $\ca_\infty^{\ell oc}$,
see the Appendix for the definition.
With certain care we shall study also the limiting situations
$\max (q_0,q_1) = \infty$.

\begin{T}\label{intera}
Let $0 < \Theta < 1$.
Let $0 < p_0,p_1 < \infty $, $0 < q_0,q_1 \le \infty$ and $s_0,s_1\in \R$. We put
\beq\label{eq:inter:1}
\frac 1 p := \frac{1-\Theta}{p_0} + \frac{\Theta}{p_1}\, , \qquad
\frac 1 q := \frac{1-\Theta}{q_0} + \frac{\Theta}{q_1}
\qquad \mbox{and}\qquad s:= (1-\Theta)\, s_0 + \Theta\, s_1\, .
\eeq
Let $w_0, w_1\in {\mathcal A}_\infty^{\ell oc}$  and define $w$ by the formula
\be\label{eq:inter:2}
w := w_{0}^{\frac{(1-\Theta)p}{p_0}}w_1^{\frac{\Theta p}{p_1}}.
\ee
Then
\be
\label{eq-28}
f^{s_0}_{p_0,q_0}(\R^d,w_0)^{1-\Theta}\,  f^{s_1}_{p_1,q_1}(\R^d,w_1)^\Theta  =
f^{s}_{p,q}(\R^d,w)
\ee
holds in the sense of equivalent quasi-norms.
\end{T}

\noindent
\bpr
By Lemma \ref{stein2} the weight $w$ belongs to $\ca_\infty^{\ell oc}$, i.e.,
our sequence spaces $f^{s}_{p,q}(\R^d,w)$ are well-defined.
Muckenhoupt weights and therefore also local Muckenhoupt weights can not vanish on a set of positive measure.
Hence, (\ref{gewicht}) holds for $w_0,w_1$ and $w$.
\\
{\em Step 1.}
We shall prove
\[
f^{s_0}_{p_0,q_0}(\R^d,w_0)^{1-\Theta} \, f^{s_1}_{p_1,q_1}(\R^d,w_1)^\Theta  \hookrightarrow
f^{s}_{p,q}(\R^d,w)\, .
\]
We suppose that  sequences $\lambda := (\lambda_{j,k})_{j,k}$, $\lambda^\ell := (\lambda_{j,k}^\ell)_{j,k}$,
$\ell =0,1$, are given and that
\[
|\lambda_{j,k}|\le |\lambda_{j,k}^0|^{1-\Theta}\cdot |\lambda^1_{j,k}|^\Theta \qquad
\]
holds for all $j\in\N_0$ and all $k\in\Z^d$.
We have to show that  there exists a constant $c$ s.t.
\[
\|\lambda|f^{s}_{p,q}(\R^d,w)\|\le c\,
\|\lambda^0|f^{s_0}_{p_0,q_0}(\R^d,w_0)\|^{1-\Theta}\cdot\|\lambda^1| f^{s_1}_{p_1,q_1}(\R^d,w_1)\|^\Theta
\]
holds for all such $\lambda, \lambda^0, \lambda^1$.
But this follows directly by H\"older's inequality (with $c=1$).
\\
{\em Step 2.} Now we turn to the proof of
\[
f^{s}_{p,q}(\R^d,w)
\hookrightarrow f^{s_0}_{p_0,q_0}(\R^d,w_0)^{1-\Theta} \, f^{s_1}_{p_1,q_1}(\R^d,w_1)^\Theta
\, .
\]
We assume in addition $\max (q_0,q_1) <\infty$.
Let the  sequence $\lambda \in f^s_{p,q}(\R^d,w)$
be given. We have to find sequences $\lambda^0$ and $\lambda^1$ such
that $|\lambda_{j,k}|\le |\lambda_{j,k}^0|^{1-\Theta}\cdot |\lambda^1_{j,k}|^\Theta$
for every $j\in\N_0$ and $k\in\Z^d$ and
\beq\label{eq:int2}
\|\lambda^0|f^{s_0}_{p_0,q_0}(\R^d,w_0)\|^{1-\Theta}\cdot\|\lambda^1| f^{s_1}_{p_1,q_1}(\R^d,w_1)\|^\Theta
\le c\,\|\lambda|f^{s}_{p,q}(\R^d,w)\|
\eeq
with some constant $c$ independent  of $\lambda$.
We follow  ideas of the proof of Theorem 8.2 in Frazier and Jawerth \cite{FJ2}, see also Bownik \cite{Bow2}.
Since $w_0,w_1$ are local Muckenhoupt weights, they are  positive and finite a.e., hence, also $w$ is positive and finite a.e..
Let
\[
A := \left\{x\in\R^d: \: 0 < \frac{w(x)}{w_0(x)}<\infty\quad \text{and}\quad 0 < \frac{w(x)}{w_1(x)}<\infty\, \right\}\, .
\]
The functions  $w$, $w_0$ and $w_1$ are locally integrable and positive a.e., therefore  $\R^d\setminus A$ is a set of  measure zero.
We  put
\[
A_\ell
:=\left\{x\in A :\left(\sum_{j=0}^\infty\sum_{k\in\Z^d} 2^{jsq}|\lambda_{j,k}|^q\cx_{j,k}(x)\right)^{1/q}
\cdot\left(\frac{w(x)}{w_0(x)}\right)^{\frac{1}{p_0}\cdot\frac{1}{\frac{p}{p_0}-\frac{q}{q_0}}}>2^\ell\right\}\, ,
\]
$ \ell\in\Z$. Obviously $A_{\ell+1} \subset A_{\ell}$, $\ell \in \Z$.
Now we introduce a (partial) decomposition of $\N_0 \times \Z^d$ by taking
\[
C_\ell := \left\{(j,k):\: |Q_{j,k}\cap A_\ell| > \frac{|Q_{j,k}|}{2}\quad\text{and}\quad
|Q_{j,k}\cap A_{\ell+1}| \le \frac{|Q_{j,k}|}{2}\right\},\quad \ell\in\Z.
\]
The sets $C_\ell$ are pairwise disjoint, i.e.,
$C_\ell \cap C_m = \emptyset$ if $\ell \neq m$.
\\
{\em Substep 2.1.} We claim that $\lambda_{j,k}=0$ holds for all tuples  $(j,k) \not \in \bigcup_\ell C_\ell$.
Let us consider one such tuple $(j_0,k_0)$ and let us choose $l_0\in \Z$ arbitrarily.
As $(j_0,k_0) \not \in  C_{\ell_0}$, then either
\be\label{nochwas1}
|Q_{j_0,k_0}\cap A_{\ell_0}| \le  \frac{|Q_{j_0,k_0}|}{2}\qquad\text{or}\qquad
|Q_{j_0,k_0}\cap A_{\ell_0+1}| > \frac{|Q_{j_0,k_0}|}{2}\, .
\ee
Let us assume for the moment that the second  condition is satisfied.
By induction on $\ell$ it follows
\be\label{2.11.2012}
|Q_{j_0,k_0}\cap A_{\ell+1}| > \frac{|Q_{j_0,k_0}|}{2}\, \qquad \mbox{for all}\quad \ell \ge \ell_0\, .
\ee
Let $D := \bigcap_\ell  Q_{j_0,k_0}\cap A_{\ell}$. The family $\{Q_{j_0,k_0}\cap A_{\ell}\}_\ell$ is a
decreasing family of sets, i.e., $Q_{j_0,k_0}\cap A_{\ell+1}\subset Q_{j_0,k_0}\cap A_{\ell}$. Therefore, in view of (\ref{2.11.2012}),
the measure of the set $D$ is larger than  or equal to   $\frac{|Q_{j_0,k_0}|}{2}$.
We obtain
\beqq
\| \, \lambda \, |f^s_{p,q}(\R^d, w)\|^p & := &
\bigg\|\bigg( \sum_{j=0}^\infty  \sum_{k\in\Z^d}
 2^{sjq}\, |\lambda_{j,k}|^q\, \cx_{j,k}(\cdot)\bigg)^{1/q}\bigg| L_p(\R^d,w)
\bigg\|^p
\\
& > & \int_{Q_{j_0,k_0}\cap A_\ell} \Big(2^\ell \Big(\frac{w_0(x)}{w(x)}\Big)^{\frac{1}{p_0}\cdot\frac{1}{\frac{p}{p_0}-\frac{q}{q_0}}}\Big)^p
w(x)\, dx
\\
& > & 2^{\ell p} \, \int_{D}  w_0(x)^{\frac{p}{p_0}\cdot\frac{1}{\frac{p}{p_0}-\frac{q}{q_0}}}\,
w(x)^{\frac{q}{q_0}\cdot\frac{1}{\frac{q}{q_0}-\frac{p}{p_0}}} \, dx\, .
\eeqq
The norm $\| \, \lambda \, |f^s_{p,q}(\R^d, w)\|^p$ is finite since $\lambda \in f^s_{p,q}(\R^d, w)$. In consequence the integral over $D$ is a finite positive number. We recall that  the function we integrate is positive a.e. and $|D|\ge \frac{|Q_{j_0,k_0}|}{2}$.
Letting $\ell$ tend to infinity we get a contradiction.
Hence, we have to turn in (\ref{nochwas1}) to the situation where the first condition is satisfied.
We claim
\[
|Q_{j_0,k_0}\cap A_{\ell}| \le  \frac{|Q_{j_0,k_0}|}{2}\qquad \mbox{for all}\qquad \ell \le \ell_0\, .
\]
Again this follows by induction on $\ell$ using $(j_0,k_0) \not \in  \bigcup_\ell C_{\ell}$.
Obviously this yields
\be\label{2.11.2012a}
|Q_{j_0,k_0}\cap A_{\ell}^c| \ge 2^{-j_0d-1}\qquad \mbox{ for all } \qquad \ell \le \ell_0 \, .
\ee
Let now $E := \bigcap_\ell  Q_{j_0,k_0}\cap A \cap A^c_{\ell}$. The family $\{Q_{j_0,k_0}\cap A\cap A^c_{\ell}\}_\ell$
satisfies
\[
(Q_{j_0,k_0}\cap A \cap  A_{\ell-1}) \quad \subset \quad (Q_{j_0,k_0}\cap A\cap A_{\ell})\, .
\]
Therefore, in view of (\ref{2.11.2012a}), the measure of the set $E$ is larger than  or equal to   $\frac{|Q_{j_0,k_0}|}{2}$.
By  selecting a point $x \in E$ we conclude that
\[
2^{j_0s}\,  |\lambda_{j_0,k_0}| \le
\left(\sum_{j=0}^\infty\sum_{k\in\Z^d} 2^{jsq}|\lambda_{j,k}|^q\cx_{j,k}(x)\right)^{1/q}
\le 2^{\ell} \, \left(\frac{w_0(x)}{w(x)}\right)^{\frac{1}{p_0}\cdot\frac{1}{\frac{p}{p_0}-\frac{q}{q_0}}}
\, ,
\]
for any $\ell\le \ell_0$. Now, for $\ell$ tending to $-\infty$ the claim, namely $\lambda_{j_0,k_0}=0$, follows.
\\
{\em Substep 2.2.}
If $(j,k)\not\in\bigcup_{\ell\in\Z}C_\ell$, then we define $\lambda_{j,k}^0 :=\lambda^1_{j,k}:=0.$
If $(j,k)\in C_\ell$, we put
\[
 \lambda_{j,k}^0 := 2^{\ell \gamma}\, 2^{ju}\, |\lambda_{j,k}|^{q/q_0}\quad\text{and}\quad
 \lambda_{j,k}^1 := 2^{\ell\delta}\, 2^{jv}\, |\lambda_{j,k}|^{q/q_1},
\]
where
\[
 \gamma :=\frac{p}{p_0}-\frac{q}{q_0}\, , \quad \delta :=\frac{p}{p_1}-\frac{q}{q_1}\, ,\quad
 u := q\, \Theta\, \left[\frac{s_1}{q_0}-\frac{s_0}{q_1}\right]\, ,
\quad
 v:= q \, (1-\Theta)\, \left[\frac{s_0}{q_1}-\frac{s_1}{q_0}\right].
\]
We observe, that
\[
(\lambda^0_{j,k})^{1-\Theta}\cdot (\lambda^1_{j,k})^\Theta=
2^{\ell[\gamma(1-\Theta)+\delta\Theta]}\cdot 2^{j[u(1-\Theta)+v\Theta]}\cdot |\lambda_{j,k}|
=|\lambda_{j,k}|\, ,
\]
which holds now for all pairs $(j,k)$.
To prove \eqref{eq:int2}, it will be sufficient to establish the following two inequalities
\beq\label{eq:int4}
 \|\lambda^0|f^{s_0}_{p_0,q_0}(\R^d,w_0)\| & \lesssim    &
 \|\lambda|f^{s}_{p,q}(\R^d,w)\|^{p/p_0}
\\
\label{eq:int4b}
 \|\lambda^1|f^{s_1}_{p_1,q_1}(\R^d,w_1)\| & \lesssim   &
 \|\lambda|f^{s}_{p,q}(\R^d,w)\|^{p/p_1}\, .
\eeq
{\em Substep 2.3}. First we deal with  \eqref{eq:int4} under the condition
$\gamma \ge 0$.
Our restrictions on $p_0,p_1$ and $q_0,q_1$ are symmetric. 
It follows from (\ref{eq:inter:1}) that
\[\min \Big(\frac{p_0}{q_0}, \, \frac{p_1}{q_1} \Big) \le  \frac pq \le  \max \Big(\frac{p_0}{q_0}, \, \frac{p_1}{q_1} \Big)\, .\]
As $\gamma\ge 0$, we get also $p/p_0  \ge q/q_0$ and $\delta \le 0$.
By employing the sets $C_\ell$ we derive
\begin{align*}
\|\lambda^0|f^{s_0}_{p_0,q_0}(\R^d,w_0)\|&=
\Big\|
\Big(\sum_{j=0}^\infty\sum_{k\in\Z^d} 2^{js_0q_0}(\lambda_{j,k}^0)^{q_0}\cx_{j,k}(\cdot)\Big)^{1/q_0}
\Big|L_{p_0}(\R^d,w_0)\Big\|\\
& = \Big\|\Big(\sum_{\ell=-\infty}^\infty \sum_{(j,k)\in C_\ell}
2^{js_0q_0} \, 2^{\ell \gamma q_0}\, 2^{juq_0}\,|\lambda_{j,k}|^q\cx_{j,k}(\cdot)
\Big)^{1/q_0} \Big|L_{p_0}(\R^d,w_0)\Big\|\\
& = \Big\|\Big(\sum_{\ell = -\infty}^\infty \sum_{(j,k)\in C_\ell}
f_{j,k} (\, \cdot \, )^Q   \Big)^{1/Q} \Big|L_{P}(\R^d,w_0)\Big\|^{P/p_0},
\end{align*}
where
\[
f_{j,k}(\, \cdot\, ) : = \Big(2^{js_0q_0} \, 2^{\ell \gamma q_0}\, 2^{juq_0}\, |\lambda_{j,k}|^q\, \cx_{j,k}(\, \cdot\, )
 \Big)^{\frac{p_0}{q_0P}},\qquad (j,k)\in C_\ell\, ,
\]
and $P$ and $Q=\frac{q_0P}{p_0}$ are chosen such that $w_0\in {\cal A}_P, 1<P<\infty, 1<Q\le\infty$.
Next we apply the weighted vector-valued maximal inequality (\ref{slava})
for the local Hardy-Littlewood maximal function  from the  Appendix together with the estimate
\[
\cx_{j,k}(x)\le 2\, (M^{\ell oc}\cx_{Q_{j,k}\cap A_\ell})(x)\, ,   \qquad  x \in \Rd, \quad (j,k)\in C_\ell\, .
\]
Using $u+s_0=\frac{sq}{q_0}$,
$\gamma\ge 0$ and
\[
\bigcup_{L=-\infty}^\infty (A_L \setminus A_{L+1}) = \bigcup_{\ell=-\infty}^\infty A_\ell
\]
we may further proceed
\begin{align*}
\|\lambda^0&|f^{s_0}_{p_0,q_0}(\Rd, w_0)\|
\lesssim \Big\|\Big(\sum_{\ell = -\infty}^\infty \sum_{(j,k)\in C_\ell}
2^{\ell \gamma q_0}\, 2^{juq_0}\, 2^{js_0q_0}\, |\lambda_{j,k}|^q\, \cx_{Q_{j,k}\cap A_\ell}(\cdot)
\Big)^{1/q_0} \Big|L_{p_0}(\R^d,w_0)\Big\|
\\
& = \Big\|\Big(\sum_{\ell = -\infty}^\infty \sum_{(j,k)\in C_\ell}
2^{\ell \gamma q_0}\,  2^{jsq}\, |\lambda_{j,k}|^q\cx_{Q_{j,k}\cap A_\ell}(\cdot)
\Big)^{1/q_0} \Big|L_{p_0}(\R^d,w_0)\Big\|
\\
& \le \Big\|\sum_{L=-\infty}^\infty \cx_{A_{L}\setminus A_{L+1}}(\cdot)\,
\Big( \sum_{\ell = -\infty}^\infty \sum_{(j,k)\in C_\ell}
2^{\ell \gamma q_0}\,  2^{jsq}\, |\lambda_{j,k}|^q \cx_{Q_{j,k}\cap A_\ell}(\cdot)
\Big)^{1/q_0} \Big| L_{p_0}(\R^d,w_0)\Big\|
\\
& \le \Big\| \sum_{L=-\infty}^\infty \cx_{A_{L}\setminus A_{L+1}}(\cdot)\,
\Big( \sum_{\ell = -\infty}^L \sum_{(j,k)\in C_\ell}
2^{\ell \gamma q_0}\, 2^{jsq}\, |\lambda_{j,k}|^q \cx_{j,k}(\cdot)
\Big)^{1/q_0} \Big|L_{p_0}(\R^d,w_0)\Big\|\, .
\end{align*}
Introducing the abbreviation $\displaystyle D_L:=\bigcup_{m\le L}C_m$ and using that $\gamma\ge 0$, we obtain
\begin{align*}
\|\lambda^0|f^{s_0}_{p_0,q_0}(\Rd, w_0)\|
& \lesssim
\Big\|\sum_{L=-\infty}^\infty \cx_{A_{L}\setminus A_{L+1}}(\cdot)\, 2^{L\gamma} \Big(\sum_{(j,k)\in D_L}
2^{jsq}\, |\lambda_{j,k}|^q\cx_{j,k}(\cdot)  \Big)^{1/q_0} \Big|L_{p_0}(\R^d,w_0)\Big\|
\\
& \le \Big\|\sum_{L=-\infty}^\infty \cx_{A_{L}\setminus A_{L+1}}(\cdot)\, 2^{L\gamma}
\Big(\sum_{j=0}^\infty \sum_{k\in\Z^d}
2^{jsq}\, |\lambda_{j,k}|^q \cx_{j,k}(\cdot)  \Big)^{1/q_0} \Big|L_{p_0}(\R^d,w_0)\Big\|\, .
\end{align*}
Let
\[
f(\cdot):=\Big(\sum_{j=0}^\infty \sum_{k\in\Z^d} 2^{jsq}\, |\lambda_{j,k}|^q \cx_{j,k}(\cdot)\Big)^{1/q}\, .
\]
We employ the definition of $A_L$  and find
\begin{align*}
\|\lambda^0|f^{s_0}_{p_0,q_0}(\R^d,w_0)\|
&\lesssim
\Big\|\sum_{L=-\infty}^\infty \cx_{A_{L}\setminus A_{L+1}}(\cdot)\, f^{\gamma}(\cdot)\,
\Big(\frac{w(\cdot)}{w_0(\cdot)}\Big)^{\frac 1{p_0}\cdot\frac{1}{\frac{p}{p_0}-\frac{q}{q_0}}\cdot\gamma} \, f^{\frac{q}{q_0}}(\cdot)
\Big| L_{p_0}(\R^d,w_0)\Big\|
\\
& = \Big\|\, f^{\gamma+q/q_0}\, \Big(\frac{w}{w_0}\Big)^{\frac 1{p_0}}
\Big| L_{p_0}(\R^d,w_0)\Big\|
\\
& = \| \, f^{\frac p{p_0}} \, (w/w_0)^{\frac 1{p_0}}\, |L_{p_0}(\R^d,w_0)\|
\\
& = \|f|L_p(\R^d,w)\|^{p/p_0} = \|\lambda|f^s_{p,q}(\R^d,w)\|^{p/p_0}.
\end{align*}
{\em Substep 2.4}. Now we prove (\ref{eq:int4b}). We only make some comments concerning necessary modifications in comparison with Substep 2.3.
We first point out, that the identity
$$
\left(\frac{w(x)}{w_0(x)}\right)^{\frac{1-\Theta}{p_0}}=
\left(\frac{w(x)}{w_1(x)}\right)^{-\frac{\Theta}{p_1}}, \qquad x\in A
$$
raised to the appropriate power gives
\beq\label{eq:int3}
\left(\frac{w(x)}{w_0(x)}\right)^{\frac{1}{p_0}\cdot\frac{1}{\frac{p}{p_0}-\frac{q}{q_0}}}
=\left(\frac{w(x)}{w_1(x)}\right)^{\frac{1}{p_1}\cdot\frac{1}{\frac{p}{p_1}-\frac{q}{q_1}}}, \qquad x\in A.
\eeq
This means, that the definition of the sets $A_\ell$ and $C_\ell$ does not change, if we replace
$(w_0, p_0,q_0)$ by $(w_1, p_1, q_1)$. As $\delta\le 0$ in this case,
we are forced to replace the sets $Q_{j,k}\cap A_\ell$ by $Q_{j,k}\cap A_{\ell+1}^c$. Observe that  $|Q_{j,k}\cap A_{\ell+1}^c| \ge \frac{|Q_{j,k}|}{2}$
and hence
\[
\cx_{j,k}(x)\le 2 \, (M^{\ell oc} \cx_{Q_{j,k}\cap A_{\ell+1}^ c})(x)\, ,   \qquad x\in \Rd, \quad  (j,k)\in C_\ell\, .
\]
This, together with $v+s_1=sq/q_1$ and the maximal inequality (\ref{slava}), leads to
\begin{align*}
\|\lambda^1&|f^{s_1}_{p_1,q_1}(\Rd, w_1)\|
\lesssim \Big\|\Big(\sum_{\ell = -\infty}^\infty \sum_{(j,k)\in C_\ell}
2^{\ell \delta q_1}\, 2^{jvq_1}\, 2^{js_1q_1}\, |\lambda_{j,k}|^q\,
\cx_{Q_{j,k}\cap A_{\ell+1}^c}(\cdot)
\Big)^{1/q_1} \Big|L_{p_1}(\R^d,w_1)\Big\|
\\
& = \Big\|\Big(\sum_{\ell = -\infty}^\infty \sum_{(j,k)\in C_\ell}
2^{\ell \delta q_1}\,  2^{jsq}\, |\lambda_{j,k}|^q
\cx_{Q_{j,k}\cap A_{\ell+1}^c}(\cdot)
\Big)^{1/q_1} \Big|L_{p_1}(\R^d,w_1)\Big\|
\\
& \le \Big\|\sum_{L=-\infty}^\infty \cx_{A_{L+1}^ c \setminus A_{L}^ c}(\cdot)\,
\Big( \sum_{\ell = -\infty}^\infty \sum_{(j,k)\in C_\ell}
2^{\ell \delta q_1}\,  2^{jsq}\, |\lambda_{j,k}|^q
\cx_{Q_{j,k}\cap A_{\ell+1}^ c}(\cdot)
\Big)^{1/q_1} \Big| L_{p_1}(\R^d,w_1)\Big\|
\\
& \le \Big\| \sum_{L=-\infty}^\infty \cx_{A_{L+1}^c\setminus A_{L}^c}(\cdot)\,
\Big( \sum_{\ell = L-1}^\infty \sum_{(j,k)\in C_\ell}
2^{\ell \delta q_1}\, 2^{jsq}\, |\lambda_{j,k}|^q
\cx_{j,k}(\cdot)\Big)^{1/q_1} \Big|L_{p_1}(\R^d,w_1)\Big\|\, .
\end{align*}
Defining $\displaystyle E_L:=\bigcup_{m\ge L-1}C_m$
and making use of $\delta \le 0$ we obtain
\[
\|\lambda^1|f^{s_1}_{p_1,q_1}(\Rd, w_1)\|
 \lesssim
\Big\|\sum_{L=-\infty}^\infty \cx_{A_{L+1}^c\setminus A_{L}^ c}(\cdot)\, 2^{L\delta} \Big(\sum_{(j,k)\in E_L}
2^{jsq}\, |\lambda_{j,k}|^q\cx_{j,k}(\cdot)  \Big)^{1/q_1} \Big|L_{p_1}(\R^d,w_1)\Big\|\, .
\]
As above we put
\[
f(\cdot):=\Big(\sum_{j=0}^\infty \sum_{k\in\Z^d} 2^{jsq}\, |\lambda_{j,k}|^q \cx_{j,k}(\cdot)\Big)^{1/q}\, .
\]
The definition of $A_L$  yields
\[
2 ^L < f(x) \, \left(\frac{w(x)}{w_1(x)}\right)^{\frac{1}{p_1}\cdot\frac{1}{\frac{p}{p_1}-\frac{q}{q_1}}}\le 2^{L+1}\, , \qquad x \in A_{L+1}^c\setminus A_{L}^ c\, ,
\]
see (\ref{eq:int3}). Inserting this in the previous inequality we get
\begin{align*}
\|\lambda^1|f^{s_1}_{p_1,q_1}(\R^d,w_1)\|
&\lesssim
\Big\|\sum_{L=-\infty}^\infty \cx_{A_{L+1}^c \setminus A_{L}^ c}(\cdot)\,
f^{\delta + \frac{q}{q_1}}(\cdot)\, \Big(\frac{w(\cdot)}{w_1(\cdot)}\Big)^{\frac 1{p_1}\cdot\frac{1}{\frac{p}{p_1}-\frac{q}{q_1}}\cdot\delta} \,
\Big| L_{p_1}(\R^d,w_1)\Big\|
\\
& \lesssim \| \, f^{\frac p{p_1}} \, (w/w_1)^{\frac 1{p_1}}\, |L_{p_1}(\R^d,w_1)\|
\\
& = \|f|L_p(\R^d,w)\|^{p/p_1} = \|\lambda|f^s_{p,q}(\R^d,w)\|^{p/p_1}.
\end{align*}
This proves (\ref{eq:int4b}).
\\
{\em Step 3.} Let $\max (q_0,q_1) = \infty$. \\
{\em Substep 3.1}.
First we consider $0 < q_1 < q_0 = \infty$.
We shall discuss
the needed modifications only. As in Step 2,
if $(j,k)\not\in\bigcup_{\ell\in\Z}C_\ell$, then we define $\lambda_{j,k}^0 =\lambda^1_{j,k}=0.$
If $(j,k)\in C_\ell$, we put
\be\label{gut4}
 \lambda_{j,k}^0 := 2^{\ell \gamma}\, 2^{ju}\quad\text{and}\quad
 \lambda_{j,k}^1 := 2^{\ell\delta}\, 2^{jv}\, |\lambda_{j,k}|^{q/q_1},
\ee
where
\[
 \gamma :=\frac{p}{p_0}\, , \quad \delta :=\frac{p}{p_1}-\frac{q}{q_1}\, ,\quad
 u :=  - q\, \Theta\,  \frac{s_0}{q_1}\, ,
\quad
 v:= q \, (1-\Theta)\, \frac{s_0}{q_1}\, .
\]
This implies
$(\lambda^0_{j,k})^{1-\Theta}\cdot (\lambda^1_{j,k})^\Theta=|\lambda_{j,k}|$.
Again we are going to establish the inequalities (\ref{eq:int4}) and (\ref{eq:int4b}).
Since there is nothing changed with respect to (\ref{eq:int4b}) we deal with  \eqref{eq:int4}.
Obviously $\gamma > 0$.
Formally it looks like that we lost the influence of $\lambda$. However, this is not true.
By employing the same arguments as in Substep 2.3 and $u+s_0 =0$ we obtain
 \begin{align*}
\|\lambda^0|f^{s_0}_{p_0,\infty}(\Rd, w_0)\|
&\lesssim \|\sup_{\ell \in \Z} \sup_{(j,k)\in C_\ell}
2^{\ell \gamma}\,  \cx_{Q_{j,k}\cap A_\ell}(\cdot)\,  |L_{p_0}(\R^d,w_0)\|
\\
& \le \Big\| \sum_{L=-\infty}^\infty \cx_{A_{L}\setminus A_{L+1}}(\cdot)\,
\sup_{-\infty < \ell \le L} \sup_{(j,k)\in C_\ell}
2^{\ell \gamma}\,  \cx_{j,k}(\cdot) \,  \Big|L_{p_0}(\R^d,w_0)\Big\|
\\
& \le
\Big\|\sum_{L=-\infty}^\infty \cx_{A_{L}\setminus A_{L+1}}(\cdot)\, 2^{L\gamma} \sup_{(j,k)\in D_L} \cx_{j,k}(\cdot)
 \Big|L_{p_0}(\R^d,w_0)\Big\|
\end{align*}
Next we use the definition of the set $A_L$.  In case $x \in \cx_{A_{L}\setminus A_{L+1}}$ we conclude
\[
2^L< \Big(\sum_{j=0}^\infty \sum_{k\in\Z^d}
2^{jsq}\, |\lambda_{j,k}|^q \, \cx_{j,k}(x)\Big)^{1/q}\, \Big(\frac{w(x)}{w_0 (x)}\Big)^{1/p} \le
2^{L+1} \, .
\]
We insert this into the previous inequality and find
\begin{align*}
\|\lambda^0&|f^{s_0}_{p_0,\infty}(\Rd, w_0)\|
\\
& \lesssim
 \Big\|\sum_{L=-\infty}^\infty \cx_{A_{L}\setminus A_{L+1}}(\cdot)\,
\Big(\sum_{j=0}^\infty \sum_{k\in\Z^d}
2^{jsq}\, |\lambda_{j,k}|^q \, \cx_{j,k}(x)\Big)^{\gamma/q}
\Big(\frac{w(x)}{w_0 (x)}\Big)^{\gamma/p}
\Big|L_{p_0}(\R^d,w_0)\Big\|
\\
& \le
\Big\|\Big(\sum_{j=0}^\infty \sum_{k\in\Z^d}
2^{jsq}\, |\lambda_{j,k}|^q \, \cx_{j,k}(x)\Big)^{\gamma/q} \Big(\frac{w(x)}{w_0 (x)}\Big)^{1/p_0} \Big|L_{p_0}(\R^d,w_0)\Big\|
\\
& =
\Big\|\Big(\sum_{j=0}^\infty \sum_{k\in\Z^d}
2^{jsq}\, |\lambda_{j,k}|^q \, \cx_{j,k}(x)\Big)^{1/q} \Big|L_{p}(\R^d,w)\Big\|^{p/p_0}\, ,
\end{align*}
as we wanted to prove. \\
{\em Substep 3.2}.
Observe that the case  $0 < q_0 < q_1 = \infty$ follows by symmetry.
\\
{\em Substep 3.3}.
It remains to study the case $q_0 = q_1 = \infty$.
If $(j,k)\in C_\ell$, we choose $\lambda_{j,k}^i$, $i=0,1$, s.t.
\be\label{gut3}
 \lambda_{j,k}^0 := 2^{\ell \gamma}\, 2^{ju} \, |\lambda_{j,k}| \quad\text{and}\quad
 \lambda_{j,k}^1 := 2^{\ell\delta}\, 2^{jv}\, |\lambda_{j,k}|,
\ee
where
\[
 \gamma :=\frac{p}{p_0} -1\, , \quad \delta :=\frac{p}{p_1}-1\, ,\quad
 u :=  \Theta\,  (s_1 - s_0)\, ,
\quad
 v:=  (1-\Theta)\, (s_0-s_1)\, .
\]
Again this implies
$(\lambda^0_{j,k})^{1-\Theta}\cdot (\lambda^1_{j,k})^\Theta=|\lambda_{j,k}|$.
Without loss of generality we may assume that $p_0 \le p \le p_1$, i.e., $\gamma \ge 0$.
Now, using $u+s_0 = s$, we may proceed as in Substep 3.1 and obtain
\[
\|\lambda^0|f^{s_0}_{p_0,\infty}(\Rd, w_0)\|
\lesssim
\Big\|\sum_{L=-\infty}^\infty \cx_{A_{L}\setminus A_{L+1}}(\cdot)\, 2^{L\gamma} \sup_{(j,k)\in D_L} 2^{js}\, |\lambda_{j,k}| \, \cx_{j,k}(\cdot)
 \Big|L_{p_0}(\R^d,w_0)\Big\|\, .
\]
Employing the definition of the set $A_L$  we conclude
\[
2^L< \sup_{j\in\N_0 }\,  \sup_{k\in\Z^d}\,
2^{js}\, |\lambda_{j,k}| \, \cx_{j,k}(x) \, \Big(\frac{w(x)}{w_0 (x)}\Big)^{\frac{1}{p_0} \, \frac{1}{\frac{p}{p_0}-1}} \le
2^{L+1} \, , \qquad x \in \cx_{A_{L}\setminus A_{L+1}}\, .
\]
We insert this into the previous inequality and find
\begin{align*}
\|\lambda^0&|f^{s_0}_{p_0,\infty}(\Rd, w_0)\|
\\
& \lesssim
 \Big\|\sum_{L=-\infty}^\infty \cx_{A_{L}\setminus A_{L+1}}(\cdot)\,
\Big(\sup_{j\in\N_0 }\,  \sup_{k\in\Z^d}
2^{js}\, |\lambda_{j,k}| \, \cx_{j,k}(x)\Big)^{\gamma + 1}
\Big(\frac{w(x)}{w_0 (x)}\Big)^{\gamma \, \frac{1}{p_0} \, \frac{1}{\frac{p}{p_0}-1}}
\Big|L_{p_0}(\R^d,w_0)\Big\|
\\
& \le
\Big\|\Big(\sup_{j\in\N_0 }\,  \sup_{k\in\Z^d}
2^{js}\, |\lambda_{j,k}| \, \cx_{j,k}(x)\Big)^{p/p_0} \Big(\frac{w(x)}{w_0 (x)}\Big)^{\frac{1}{p_0}}\,  \Big|L_{p_0}(\R^d,w_0)\Big\|
\\
& =
\Big\|\, \sup_{j\in\N_0}\,  \sup_{k\in\Z^d}
2^{js}\, |\lambda_{j,k}| \, \cx_{j,k}(x)\,  \Big|L_{p}(\R^d,w)\Big\|^{p/p_0}\, ,
\end{align*}
as we wanted to prove.
The estimate of $\|\lambda^1|f^{s_1}_{p_1,\infty}(\Rd, w_1)\|$ follows by similar arguments.
\epr

\begin{Rem}
 \rm
(i) The identity
\[
f^{s_0}_{p_0,q_0}(\R^d)^{1-\Theta}\,  f^{s_1}_{p_1,q_1}(\R^d)^\Theta  =
f^{s}_{p,q}(\R^d)
\]
i.e., formula (\ref{eq-28}) with $w_0 = w_1 \equiv 1$, has been proved in Frazier and Jawerth \cite{FJ2}.
Our proof, given here, is just an adaptation to the weighted situation.
However, let us mention, that we had some advantage from Bownik's proof in \cite{Bow2}, in particular
in case $\max (q_0,q_1) = \infty$.
In fact, Bownik had considered the  situation
\[
f^{s_0}_{p_0,q_0}(\R^d, \mu)^{1-\Theta}\,  f^{s_1}_{p_1,q_1}(\R^d,\mu)^\Theta  =
f^{s}_{p,q}(\R^d,\mu)\, ,
\]
where $\mu$ is a doubling measure. In Substep 2.1 we also used some ideas from Yang, Yuan and Zhuo \cite{DWC}.
These three authors dealt with extensions to
sequence spaces related to Lizorkin-Triebel spaces built on Morrey spaces.
Further we would like to mention that
Wojciechowska \cite{AW2} recently  proved
\[
f^{s_0}_{p_0,q_0}(\R^d, w)^{1-\Theta} \, f^{s_1}_{p_1,q_1}(\R^d,w)^\Theta =
f^{s}_{p,q}(\R^d,w)\, ,
\]
where $w \in \ca_\infty^{\ell oc}$. The class $\ca_\infty^{\ell oc}$
and the set of doubling measures are incomparable
(more exactly, there exists weights in  $\ca_\infty^{\ell oc}$ which are exponentially
growing and therefore do not induce a doubling measure, vice versa there are doubling measures
which do not induce a weight in   $\ca_\infty$). An example of a doubling measure such that the associated weight  does
not belong to $\ca_\infty^{\ell oc}$ can be found in Wik's paper \cite{Wik}. Furthermore, we refer to \cite[\S I.8.8]{stein} for an example
of a doubling measure, which is singular with respect to the Lebesgue measure, i.e. without any associated weight.\\
Theorem  \ref{intera} with $w_0 \neq w_1 $ seems to be a novelty.
However, as the previously mentioned results of Bownik  indicate, there is some hope to extend it (as well as Corollary \ref{inter})
to  larger classes of weights.
\\
(ii) Frazier, Jawerth \cite{FJ2} and Bownik \cite{Bow2} have also treated
extensions of Theorem  \ref{intera} to $\max (p_0,p_1) = \infty$.
\end{Rem}

At the end of this subsection we would like to discuss the class of admissible weights.
We had concentrated on the class $\ca^{\ell oc}_\infty$ in Theorem \ref{intera}.
We do not believe that this is the end of the story and expect that Theorem \ref{intera} holds also for more general weights.
Let us assume for the moment that Proposition \ref{sawano} extends to some weights $w_0,w_1$
(not necessarily belonging to $\ca_\infty^{\ell oc}$).
In addition we need the identification of
$L_{p_j} (\Rd,w_j)$ and $F^0_{p_j,2} (\Rd,w_j)$, $1 \le p_j <\infty$, $j=0,1$, which is known to be true
only for the class $\ca_\infty^{\ell oc}$, see Rychkov \cite{Ry}.
These two properties are also needed for the associated weight $w$.
Then Lemma \ref{beispiel} implies that
\[
F^0_{p_0,2} (\Rd,w_0)^{1-\Theta} \, F^0_{p_1,2} (\Rd,w_1)^{\Theta} = F^0_{p,2} (\Rd,w)
\]
and therefore
\[
f^{d/2}_{p_0,2} (\Rd,w_0)^{1-\Theta} \, f^{d/2}_{p_1,2} (\Rd,w_1)^{\Theta} = f^{d/2}_{p,2} (\Rd,w)\, .
\]


\subsection{Calder{\'o}n products of $b^s_{p,q}(\Rd,w)$ spaces}
\label{Calderonc}


\begin{Def}
Let $0<p,q\le \infty$, $s\in \R$ and let  $w:~ \Rd \to [0,\infty)$ be a  nonnegative and locally integrable
function.
We put
\begin{multline}
b^s_{p,q}(\R^d, w)  :=   \Bigg\{
\{\lambda_{j,k}\}_{j,k} : \quad\lambda_{j,k} \in \C\, ,
\\
\| \, (\lambda_{j,k}) \, |b^s_{p,q}(\R^d, w)\| :=
\bigg( \sum_{j=0}^\infty  \bigg\|\sum_{k\in\Z^d} 2^{sj} \, |\lambda_{j,k}|\, \cx_{j,k}(\cdot)\, \bigg| L_p(\R^d,w)
\bigg\|^q\bigg)^{1/q} \!\!\!\!< \infty \Bigg\} \, .
\end{multline}
\end{Def}

\begin{Rem}
\rm
(i) In case $w(x) = 1$ for all $x \in \Rd$ we are back in the unweighted situation.
The associated sequence spaces are denoted simply by
$b^s_{p,q}(\R^d)$.
\\
(ii)
Let $w$ satisfy (\ref{gewicht}).
Let $\mathring{b}^s_{p,q} (\Rd,w)$ denote the closure of the finite sequences in ${b}^s_{p,q} (\Rd,w)$.
We have
\[
\mathring{b}^s_{p,q} (\Rd,w) = {b}^s_{p,q} (\Rd,w)\qquad \Longleftrightarrow \qquad \max (p,q) < \infty\, .
\]
If $\max (p,q) = \infty$, then $\mathring{b}^s_{p,q} (\Rd,w)$ is a proper subspace of
${b}^s_{p,q} (\Rd,w)$. \\
(iii) Let $w$ satisfy (\ref{gewicht}).
It is easily checked that $b^s_{p,q}(\Rd,w)$ is separable if, and only if,
$\max (p,q) < \infty$.
\end{Rem}

Before we are turning to a description
of the associated Calder{\'o}n products we would like to introduce a second type of sequence spaces.
Observe
\[
\bigg\|\sum_{k\in\Z^d} 2^{js} \, |\lambda_{j,k}|\, \cx_{j,k}(\cdot)\, \bigg| L_p(\R^d,w)
\bigg\| = \Big( \sum_{k \in \Zd} 2^{jsp} \, |\lambda_{j,k}|^p \int_{Q_{j,k}} w(x)\, dx
\Big)^{1/p}\, .
\]
Now we shall replace $\int_{Q_{j,k}} w(x)\, dx$ by the positive real number $y_{j,k}$,
i.e., instead of a weight function we are using a sequence of positive real numbers.

\begin{Def}
Let $0< p, q\le \infty$, $s\in \R$ and let  $ y := (y_{j,k})_{j,k}$ be a sequence of positive real numbers.
We put
\begin{multline}
b^s_{p,q}(\R^d, s-y)  :=   \Bigg\{
\{\lambda_{j,k}\}_{j,k} : \quad\lambda_{j,k} \in \C\, ,
\\
\| \, (\lambda_{j,k}) \, |b^s_{p,q}(\R^d,s-y)\| :=
\bigg( \sum_{j=0}^\infty  \Big(\sum_{k\in\Z^d} 2^{jsp} \, |\lambda_{j,k}|^p \, y_{j,k}\, \Big)^{q/p} \bigg)^{1/q} \!\!\!\!< \infty \Bigg\} \, .
\end{multline}
\end{Def}

\begin{Rem}
\rm
Each space $b^s_{p,q} (\Rd,w)$, $w \in {\mathcal A}^{\ell oc}_\infty$, can be interpreted as a space $b^s_{p,q} (\Rd,s-y)$ by taking
\[
y_{j,k} := \int_{Q_{j,k}} w(x)\, dx\, , \qquad j \in \N_0\, , \quad k \in \Zd.
\]
\end{Rem}

As a consequence of the formula
\[
b^{s_0}_{p_0,q_0}(\R^d)^{1-\Theta} \, b^{s_1}_{p_1,q_1}(\R^d)^\Theta  =
b^{s}_{p,q}(\R^d)\, ,
\]
due to Kalton, Mayboroda, Mitrea \cite[Proposition~9.3]{KMM}, we derive the following.

\begin{Cor}\label{interc}
Let $0 < \Theta < 1$.
Let $0 < p_0,p_1 \le \infty $, $0 < q_0,q_1 \le \infty$ and $s_0,s_1\in \R$.
Let $p,q$ and $s$ be as in (\ref{eq:inter:1}).
Let $y^0:= (y^0_{j,k})_{j,k}, y^1:= (y^1_{j,k})_{j,k}$  be sequences of positive real numbers. We put
\[
y_{j,k} :=
(y_{j,k}^{0})^{\frac{(1-\Theta)p}{p_0}} \, (y_{j,k}^1)^{\frac{\Theta p}{p_1}}.
\]
Then
\be
\label{eq-27c}
b^{s_0}_{p_0,q_0}(\R^d,s-y^0)^{1-\Theta} \, b^{s_1}_{p_1,q_1}(\R^d,s-y^1)^\Theta  =
b^{s}_{p,q}(\R^d,s-y)
\ee
holds in the sense of equivalent quasi-norms.
\end{Cor}

\noindent
\bpr
Of course, $y:= (y_{j,k})_{j,k}$ is a sequence of positive real numbers.
\\
{\em Step 1}. A preparation.
Let $\varrho := (\varrho_{j,k})_{j,k}$ be a sequence of positive real numbers.
We introduce the associated  family of  mappings
\[
T_{\varrho, p} : ~ (\lambda_{j,k})_{j,k} \mapsto  (\lambda_{j,k} \, \cdot \, \varrho_{j,k}^{1/p} )_{j,k}\, ,
\]
where $0 < p \le \infty$ is fixed. Obviously, $T_{\varrho ,p}$ is an isomorphism
considered as a mapping
of $b^ s_{p,q}(\Rd,s-\varrho)$ onto $b^ s_{p,q}(\Rd)$ for all $s$ and all $q$.
\\
{\em Step 2}.
Again the embedding
\[
b^{s_0}_{p_0,q_0}(\R^d,s-y^0)^{1-\Theta} \, b^{s_1}_{p_1,q_1}(\R^d,s-y^1)^\Theta  \hookrightarrow
b^{s}_{p,q}(\R^d,s-y)
\]
follows  by repeated use of H\"older's inequality.
\\
{\em Step 3.}
We deal with
\[
b^{s}_{p,q}(\R^d,s-y)
\hookrightarrow b^{s_0}_{p_0,q_0}(\R^d,s-y^0)^{1-\Theta} \, b^{s_1}_{p_1,q_1}(\R^d,s-y^1)^\Theta
\, .
\]
Let the  sequence $\lambda \in b^s_{p,q}(\R^d,s-y)$
be given. We have to find sequences $\lambda^0$ and $\lambda^1$ such
that $|\lambda_{j,k}|\le |\lambda_{j,k}^0|^{1-\Theta}\cdot |\lambda^1_{j,k}|^\Theta$
for every $j\in\N_0$ and $k\in\Z^d$ and
\[
\|\lambda^0|b^{s_0}_{p_0,q_0}(\R^d,w_0)\|^{1-\Theta}\cdot\|\lambda^1| b^{s_1}_{p_1,q_1}(\R^d,w_1)\|^\Theta
\le c\,\|\lambda|b^{s}_{p,q}(\R^d,w)\|
\]
with some constant $c$ independent  of $\lambda$.
Now we are going to use the unweighted case
\[
b^{s}_{p,q}(\R^d)=  b^{s_0}_{p_0,q_0}(\R^d)^{1-\Theta} \, b^{s_1}_{p_1,q_1}(\R^d)^\Theta\, ,
\]
see Kalton, Mayboroda, Mitrea \cite[Proposition~9.3]{KMM}.
Let $\gamma_{j,k} := \lambda_{j,k} \, \cdot \, y_{j,k}^{1/p}$, $j \in \N_0$, $k\in \Zd$.
Since $\gamma :=(\gamma_{j,k})_{j,k} \in b^s_{p,q}(\R^d)$ this implies the existence
of sequences
$\gamma^0  :=(\gamma_{j,k}^0)_{j,k} \in b^{s_0}_{p_0,q_0}(\R^d)$,
$\gamma^1 :=(\gamma_{j,k}^1)_{j,k} \in b^{s_1}_{p_1,q_1}(\R^d)$
s.t.
$|\gamma_{j,k}|\le |\gamma_{j,k}^0|^{1-\Theta}\cdot |\gamma^1_{j,k}|^\Theta$
for all $j\in\N_0$ and all $k\in\Z^d$ and
\[
\|\gamma^0|b^{s_0}_{p_0,q_0}(\R^d)\|^{1-\Theta}\cdot
\|\gamma^1| b^{s_1}_{p_1,q_1}(\R^d)\|^\Theta
\le c\,\|\gamma|b^{s}_{p,q}(\R^d)\|
\]
with some constant $c$ independent  of $\gamma$.
We define
\[
\lambda_{j,k}^0 := \frac{\gamma_{j,k}^0}{(y^0_{j,k})^{1/p_0}}
\qquad \mbox{and}\qquad \lambda_{j,k}^1 := \frac{\gamma_{j,k}^1}{(y^1_{j,k})^{1/p_1}} \, .
\]
The sequences $\lambda^0$, $\lambda^1$ obviously meet the above requirements.
\epr

For  nonnegative and locally integrable functions $w_0$, $w_1$, satisfying (\ref{gewicht}),  we define
\[
y_{j,k}^i := \int_{Q_{j,k}} w_i(x)\, dx\, , \qquad j \in \N_0\, , \quad k \in \Zd\, , \quad i = 0, 1 \, .
\]
Then
\beqq
b^{s_0}_{p_0,q_0}(\R^d,w_0)^{1-\Theta} \, b^{s_1}_{p_1,q_1}(\R^d,w_1)^\Theta  & = &
b^{s_0}_{p_0,q_0}(\R^d,s-y^0)^{1-\Theta} \, b^{s_1}_{p_1,q_1}(\R^d,s-y^1)^\Theta
\\
& = &
b^{s}_{p,q}(\R^d,s-y)
\eeqq
follows by using the natural interpretation, see (\ref{eq-27c}).
However, without extra conditions on the functions $w_0$ and  $w_1$ we can not interpret
$b^{s}_{p,q}(\R^d,s-y)$ as the space
$b^{s}_{p,q}(\R^d,w)$, where the function $w$ is defined as in (\ref{eq:inter:2}).
A sufficient condition consists in
\[
\int_{Q_{j,k}} w (x)\, dx  \asymp \Big(\int_{Q_{j,k}} w_{0} (x) \, dx \Big)^{\frac{(1-\Theta)p}{p_0}}
\, \Big( \int_{Q_{j,k}} w_1 (x)\, dx \Big)^{\frac{\Theta p}{p_1}}
\]
for all $j$ and all $k$. This gives rise to the following definition.

\begin{Def}
Let $0 < p_0, p_1 \le \infty$ and $0 <\Theta<1 $.
Define $p$ by $1/p:= (1-\Theta)/p_0 + \Theta/p_1$.
Let   $w_j:~ \Rd \to [0,\infty)$, $j=0,1$, be   nonnegative
and locally integrable  functions s.t. (\ref{gewicht}) is satisfied.
We say that
the pair $(w_0,w_1)$ belongs to the class ${\mathfrak W}(\Theta, p_0,p_1)$
if
\be
\int_{Q_{j,k}} w_{0}(x)^{\frac{(1-\Theta)p}{p_0}} \, w_1(x)^{\frac{\Theta p}{p_1}} \, dx  \asymp \Big(\int_{Q_{j,k}}
w_{0} (x) \, dx \Big)^{\frac{(1-\Theta)p}{p_0}}
\, \Big( \int_{Q_{j,k}} w_1 (x)\, dx \Big)^{\frac{\Theta p}{p_1}}
\ee
holds for all $j \in \N_0$ and all $k\in \Zd$.
\end{Def}

\begin{Lem}\label{interq}
Suppose $w_0,w_1 \in \ca_\infty^{\ell oc}$.
Then the pair $(w_0,w_1)$ belongs to all classes ${\mathfrak W}(\Theta, p_0,p_1)$.
\end{Lem}

\noindent
\bpr
The inequality
\[
\int_{Q_{j,k}} w (x)\, dx  \le \Big(\int_{Q_{j,k}} w_{0} (x) \, dx \Big)^{\frac{(1-\Theta)p}{p_0}}
\, \Big( \int_{Q_{j,k}} w_1 (x)\, dx \Big)^{\frac{\Theta p}{p_1}}
\]
is a consequence of H\"older's inequality.
To prove the opposite inequality we shall  use the following characterization of ${\mathcal A}_\infty^{\ell oc}$:
a weight $v$ belongs to ${\mathcal A}_\infty^{\ell oc}$ if, and only if, there exists a constant $C>0$ s.t.
\[
\frac{1}{|Q|} \int_Q v(x) dx\, \le\,C\, \exp\Big( \frac{1}{|Q|} \int_Q \log v(x) dx\Big)
\]
holds for all cubes (with sides parallel to the axes) and volume $\le 1$,
cf. Theorem 2.15, p. 407,  in \cite{GC}. Using this inequality we obtain, for all such cubes $Q$,
\[
\Big(\frac{1}{|Q|} \int_Q w_0(x) dx\Big)^{\frac{(1-\Theta)p}{p_0}}\, \le\,  C_0 \, \exp \Big( \frac{(1-\Theta)p}{p_0} \, \frac{1}{|Q|} \,
\int_Q \log w_0(x)\,  dx\Big)
\]
and
\[
\Big(\frac{1}{|Q|} \int_Q w_1(x) dx\Big)^{\frac{\Theta p}{p_1}}\, \le\,  C_1 \,
\exp\Big(\frac{\Theta p}{p_1}\,  \frac{1}{|Q|} \, \int_Q \log w_1(x)\, dx\Big)\, .
\]
Multiplying these inequalities and applying Jensen's inequality  yields
\beqq
 & & \hspace*{-1.7cm}
\Big(\int_{Q} w_{0} (x) \, dx \Big)^{\frac{(1-\Theta)p}{p_0}}\, \Big( \int_{Q} w_1 (x)\, dx \Big)^{\frac{\Theta p}{p_1}}
\\
& \le  & C_0\, C_1 \,  |Q|  \exp\Big( \frac{1}{|Q|} \int_Q \log (w_0(x))^{\frac{(1-\Theta)p}{p_0}} dx +
\frac{1}{|Q|} \int_Q \log (w_1(x))^{\frac{\Theta p}{p_1}} dx\Big)
\\
& = & C\, |Q|  \exp\Big( \frac{1}{|Q|} \int_Q \log\big( (w_0(x))^{\frac{(1-\Theta)p}{p_0}} (w_1(x))^{\frac{\Theta p}{p_1}}\big)
\, dx\Big)
\\
& \le  & C\,  \int_Q (w_0(x))^{\frac{(1-\Theta)p}{p_0}} (w_1(x))^{\frac{\Theta p}{p_1}} \, dx
\\
&  = &  \, C \int_{Q} w (x)\, dx  \, ,
\eeqq
which completes the proof of the equivalence.
\epr

\begin{Cor}\label{interu}
Let $0 < \Theta < 1$.
Let $0 < p_0,p_1 \le \infty $, $0 < q_0,q_1 \le \infty$ and $s_0,s_1\in \R$.
Let $p,q$ and $s$ be as in (\ref{eq:inter:1}). Let $w_0, w_1 \in \ca_\infty^{\ell oc}$ and $w$ as in (\ref{eq:inter:2}).
Then
\be
\label{eq-27d}
b^{s_0}_{p_0,q_0}(\R^d,w_0)^{1-\Theta} \, b^{s_1}_{p_1,q_1}(\R^d,w_1)^\Theta  =
b^{s}_{p,q}(\R^d,w)
\ee
holds in the sense of equivalent quasi-norms.
\end{Cor}

\noindent
\bpr
This is an immediate consequence of Lemma \ref{interq} and Corollary \ref{interc}.
\epr

\begin{Rem}
\rm
(i) Any extension of Lemma \ref{interq} yields an extension of Corollary \ref{interu}. 
From this point of view it would be of interest to characterize the classes ${\mathfrak W}(\Theta, p_0,p_1)$.
\\
(ii) The formula
\[
b^{s_0}_{p_0,q_0}(\R^d)^{1-\Theta}\,  b^{s_1}_{p_1,q_1}(\R^d)^\Theta  =
b^{s}_{p,q}(\R^d)
\]
has been proved by Mendez and Mitrea \cite{MM} under the additional restriction
$s_0 \neq s_1$. Later on this restriction has been removed by
Kalton, Mayboroda and Mitrea \cite{KMM}.
Also in the situation of the $b$-spaces  the case of different weights seems to be new.
\end{Rem}


\subsection{Calder{\'o}n products of $\mathring{a}^s_{p,q}(\Rd,w)$ spaces}
\label{Calderone}


It is of certain use to study
Calder{\'o}n products of the spaces $\mathring{f}^s_{p,q}(\Rd,w)$
and $\mathring{b}^s_{p,q}(\Rd,w)$ separately.

\begin{T}\label{interx}
Let $0 < \Theta < 1$.
Let $0 < p_0,p_1 < \infty $, $0 < q_0,q_1 \le \infty$ and $s_0,s_1\in \R$.
Let $p,q$ and $s$ be defined as in (\ref{eq:inter:1}).
Let $w_0, w_1\in {\mathcal A}_\infty^{\ell oc}$  and define $w$ by the formula
(\ref{eq:inter:2}).
Then
\[
\mathring{f}^{s_0}_{p_0,q_0}(\R^d,w_0)^{1-\Theta}\,  \mathring{f}^{s_1}_{p_1,q_1}(\R^d,w_1)^\Theta  =
\mathring{f}^{s}_{p,q}(\R^d,w)
\]
and
\[
\mathring{f}^{s_0}_{p_0,q_0}(\R^d,w_0)^{1-\Theta}\,  {f}^{s_1}_{p_1,q_1}(\R^d,w_1)^\Theta
 = {f}^{s_0}_{p_0,q_0}(\R^d,w_0)^{1-\Theta}\,  \mathring{f}^{s_1}_{p_1,q_1}(\R^d,w_1)^\Theta
= \mathring{f}^{s}_{p,q}(\R^d,w)
\]
hold in the sense of equivalent quasi-norms.
\end{T}

\noindent
\bpr
The cases $\max (q_0,q_1)<\infty$ are already covered by Theorem \ref{intera}.
We may concentrate on $\max (q_0,q_1)=\infty$.
It will be enough to make some comments to the needed modifications
in the proof of  Theorem \ref{intera}.
\\
{\em Step 1.}
We shall prove
\[
\mathring{f}^{s_0}_{p_0,\infty}(\R^d,w_0)^{1-\Theta} \, f^{s_1}_{p_1,q_1}(\R^d,w_1)^\Theta  \hookrightarrow \mathring{f}^{s}_{p,q}(\R^d,w)\, .
\]
We suppose that  sequences $\lambda := (\lambda_{j,k})_{j,k}$, $\lambda^\ell := (\lambda_{j,k}^\ell)_{j,k}$,
$\ell =0,1$, are given and that
\begin{equation}\label{ls_23_11}
|\lambda_{j,k}|\le |\lambda_{j,k}^0|^{1-\Theta}\cdot |\lambda^1_{j,k}|^\Theta \qquad
\end{equation}
holds for all $j\in\N_0$ and all $k\in\Z^d$.
Now the essential observation is that if $\lambda^0 := (\lambda_{j,k}^0)_{j,k}$
is a finite sequence (i.e. only a finite number of the $\lambda_{j,k}^0$ is not vanishing),
then $\lambda$ has the same property.
Employing Step 1 of the proof of Theorem \ref{intera} we know that
\[
\|\lambda|f^{s}_{p,q}(\R^d,w)\|\le
\|\lambda^0|f^{s_0}_{p_0,\infty}(\R^d,w_0)\|^{1-\Theta}\cdot\|\lambda^1| f^{s_1}_{p_1,q_1}(\R^d,w_1)\|^\Theta
\]
holds for all such $\lambda, \lambda^0, \lambda^1$.

For $\lambda \in {f}^{s}_{p,q}(\R^d,w)$ we define the $M$-cutoff  sequence $\lambda^{(M)}$, $M\in \N$, by putting $\lambda^{(M)}_{j,k}=0$ if $j>M$ or $\sup_i|k_i|>M$ and $\lambda^{(M)}_{j.k}= \lambda_{j,k}$ otherwise. Then $\lambda \in \mathring{f}^{s}_{p,q}(\R^d,w)$ if and only if  $\lambda^{(M)}$ converge to $\lambda$ in   
${f}^{s}_{p,q}(\R^d,w)$ if $M\rightarrow \infty$,  cf. \cite{DWC}.  
Now \eqref{ls_23_11} implies that
\[
|\lambda_{j,k}-\lambda_{j,k}^{(M)}|\le |\lambda_{j,k}^{0} - \lambda_{j,k}^{0\, (M)}|^{1-\Theta}\cdot |\lambda^1_{j,k}|^\Theta \, ,
\]
thus 
\[
\|\, \lambda - \lambda^{(M)}\, |f^{s}_{p,q}(\R^d,w)\|\le
\|\, \lambda^0-\lambda^{0\,(M)}\, |f^{s_0}_{p_0,q_0}(\R^d,w_0)\|^{1-\Theta}\cdot\|\lambda^1| f^{s_1}_{p_1,q_1}(\R^d,w_1)\|^\Theta\, .
\]
Hence, $\lambda \in \mathring{f}^{s}_{p,q}(\R^d,w)$ and from this the claim follows.
The embedding
\[
{f}^{s_0}_{p_0,q_0}(\R^d,w_0)^{1-\Theta}\,  \mathring{f}^{s_1}_{p_1,q_1}(\R^d,w_1)^\Theta \hookrightarrow
\mathring{f}^{s}_{p,q}(\R^d,w)
\]
follows by symmetry. Finally, observe
\[
\mathring{f}^{s_0}_{p_0,q_0}(\R^d,w_0)^{1-\Theta}\,  \mathring{f}^{s_1}_{p_1,q_1}(\R^d,w_1)^\Theta
\hookrightarrow
\mathring{f}^{s}_{p,q}(\R^d,w)\, . 
\]
{\em Step 2.} Now we turn to the proof of
\[
\mathring{f}^{s}_{p,q}(\R^d,w)
\hookrightarrow \mathring{f}^{s_0}_{p_0,q_0}(\R^d,w_0)^{1-\Theta} \, \mathring{f}^{s_1}_{p_1,q_1}(\R^d,w_1)^\Theta
\, .
\]
Let the  sequence $\lambda \in \mathring{f}^s_{p,q}(\R^d,w)$
be given. We have to find sequences $\lambda^0 \in \mathring{f}^{s_0}_{p_0,q_0}(\R^d,w_0)$ and
$\lambda^1 \in \mathring{f}^{s_1}_{p_1,q_1}(\R^d,w_1)$ such
that $|\lambda_{j,k}|\le |\lambda_{j,k}^0|^{1-\Theta}\cdot |\lambda^1_{j,k}|^\Theta$
for every $j\in\N_0$ and $k\in\Z^d$ and
\be\label{ring1}
\|\lambda^0|f^{s_0}_{p_0,q_0}(\R^d,w_0)\|^{1-\Theta}\cdot\|\lambda^1| f^{s_1}_{p_1,q_1}(\R^d,w_1)\|^\Theta
\le c\,\|\lambda|f^{s}_{p,q}(\R^d,w)\|
\ee
with some constant $c$ independent  of $\lambda$.
For the moment we suppose that $\lambda$ is a finite sequence.
Since $\max (q_0,q_1)=\infty$, we have to use the formulas (\ref{gut4}) and (\ref{gut3})
from Step 3 of the proof of Theorem \ref{intera}.
In case  $q_0=q_1=\infty$ it is immediate that $\lambda^0$ and $\lambda^1$
are also finite sequences.
For the case $0 < q_1 < q_0=\infty$ we need a simple modification
(and similarly in case $0 < q_0 < q_1=\infty$). 
 Obviously we may assume that $\lambda^1$ is a finite sequence.
We define
\[
\lambda^0_{j,k} := \left\{
\begin{array}{lll}
 2^{\ell \gamma}\, 2^{ju} & \qquad & \mbox{if} \quad \lambda^1_{j,k} \neq 0,\\
0 & \qquad & \mbox{otherwise}\, .
\end{array}
\right.
\]
Then also $\lambda^0 $ is a finite sequence and (\ref{ring1}) is guaranteed by Step 3 of the proof of Theorem \ref{intera}.
Now we turn to a Cauchy sequence $(\lambda^{(M)})_M$ of finite sequences, convergence in
${f}^{s}_{p,q}(\R^d,w)$ with limit $\lambda$.
Using the formulas (\ref{gut4}) and (\ref{gut3}) we obtain a sequence $\lambda^{1(M)}$ of finite sequences
which is convergent in ${f}^{s_1}_{p_1,q_1}(\R^d,w_1)$ with limit $\lambda^1$. Similarly, taking into account the above modification,
we get a sequence $\lambda^{0(M)}$ of finite sequences
which is convergent in ${f}^{s_0}_{p_0,q_0}(\R^d,w_0)$ with limit $\lambda^0$.
In both cases convergence is derived by using the inequalities
\beqq
 \|\lambda^0|f^{s_0}_{p_0,q_0}(\R^d,w_0)\| & \lesssim    &
 \|\lambda|f^{s}_{p,q}(\R^d,w)\|^{p/p_0}
\\
 \|\lambda^1|f^{s_1}_{p_1,q_1}(\R^d,w_1)\| & \lesssim   &
 \|\lambda|f^{s}_{p,q}(\R^d,w)\|^{p/p_1}\, ,
\eeqq
see again Step 3 of the proof of Theorem \ref{intera}.
Since
\[
\mathring{f}^{s_0}_{p_0,q_0}(\R^d,w_0)^{1-\Theta} \, \mathring{f}^{s_1}_{p_1,q_1}(\R^d,w_1)^\Theta
\hookrightarrow
{f}^{s_0}_{p_0,q_0}(\R^d,w_0)^{1-\Theta} \, \mathring{f}^{s_1}_{p_1,q_1}(\R^d,w_1)^\Theta
\]
as well as
\[
\mathring{f}^{s_0}_{p_0,q_0}(\R^d,w_0)^{1-\Theta} \, \mathring{f}^{s_1}_{p_1,q_1}(\R^d,w_1)^\Theta
\hookrightarrow
\mathring{f}^{s_0}_{p_0,q_0}(\R^d,w_0)^{1-\Theta} \, {f}^{s_1}_{p_1,q_1}(\R^d,w_1)^\Theta
\]
the proof is complete.
\epr

\begin{Rem}
 \rm
In case $w=w_0=w_1\equiv 1$ we refer to Yang, Yuan and Zhuo \cite{DWC} for a partially different proof.
\end{Rem}

To derive the counterpart for the $b$-spaces we will not use the method of
Kalton, Mayboroda and Mitrea \cite{KMM}.
These authors reduced the proof of
\[
{b}^{s}_{p,q}(\R^d) = {b}^{s_0}_{p_0,q_0}(\R^d)^{1-\Theta}\,  {b}^{s_1}_{p_1,q_1}(\R^d)^\Theta
\]
to the complex interpolation formula
\[
{B}^{s}_{p,q}(\R^d) = [{B}^{s_0}_{p_0,q_0}(\R^d),\,  {B}^{s_1}_{p_1,q_1}(\R^d)]_\Theta\, .
\]
This time our aim will consists in proving complex interpolation formulas based on assertions on
Calder{\'o}n products.
Our proof will be based on the results in the unweighted case, i.e.,
\[
\mathring{b}^{s_0}_{p_0,q_0}(\R^d)^{1-\Theta}\,  \mathring{b}^{s_1}_{p_1,q_1}(\R^d)^\Theta  =
\mathring{b}^{s}_{p,q}(\R^d)
\]
and
\[
\mathring{b}^{s_0}_{p_0,q_0}(\R^d)^{1-\Theta}\,  {b}^{s_1}_{p_1,q_1}(\R^d)^\Theta
 = {b}^{s_0}_{p_0,q_0}(\R^d)^{1-\Theta}\,  \mathring{b}^{s_1}_{p_1,q_1}(\R^d)^\Theta
= \mathring{b}^{s}_{p,q}(\R^d)
\]
in the sense of equivalent quasi-norms.
For this we refer to the recent paper of Yang, Yuan and Zhuo \cite{DWC}.\\

\begin{T}\label{intery}
Let $0 < \Theta < 1$.
Let $0 < p_0,p_1, q_0,q_1 \le \infty$ and $s_0,s_1\in \R$.
Let $p,q$ and $s$ be defined as in (\ref{eq:inter:1}).
Let $w_0, w_1\in {\mathcal A}_\infty^{\ell oc}$  and define $w$ by the formula
(\ref{eq:inter:2}).
Then
\be
\label{eq-28b}
\mathring{b}^{s_0}_{p_0,q_0}(\R^d,w_0)^{1-\Theta}\,  \mathring{b}^{s_1}_{p_1,q_1}(\R^d,w_1)^\Theta  =
\mathring{b}^{s}_{p,q}(\R^d,w)
\ee
and
\be
\label{eq-28c}
\mathring{b}^{s_0}_{p_0,q_0}(\R^d,w_0)^{1-\Theta}\,  {b}^{s_1}_{p_1,q_1}(\R^d,w_1)^\Theta
 = {b}^{s_0}_{p_0,q_0}(\R^d,w_0)^{1-\Theta}\,  \mathring{b}^{s_1}_{p_1,q_1}(\R^d,w_1)^\Theta
= \mathring{b}^{s}_{p,q}(\R^d,w)
\ee
hold in the sense of equivalent quasi-norms.
\end{T}

\noindent
\bpr
{\em Step 1.} Concerning  the embedding
\[
\mathring{b}^{s_0}_{p_0,q_0}(\R^d,w_0)^{1-\Theta} \, b^{s_1}_{p_1,q_1}(\R^d,w_1)^\Theta  \hookrightarrow
\mathring{b}^{s}_{p,q}(\R^d,w)
\]
the arguments from Step 1 of the proof of Theorem \ref{interx} carry over to the present situation.
The embedding
\[
{b}^{s_0}_{p_0,q_0}(\R^d,w_0)^{1-\Theta} \, \mathring{b}^{s_1}_{p_1,q_1}(\R^d,w_1)^\Theta  \hookrightarrow
\mathring{b}^{s}_{p,q}(\R^d,w)
\]
follows by symmetry.
\\
{\em Step 2.}
It remains to prove
\[
\mathring{b}^{s}_{p,q}(\R^d,w)
\hookrightarrow \mathring{b}^{s_0}_{p_0,q_0}(\R^d,w_0)^{1-\Theta} \, \mathring{b}^{s_1}_{p_1,q_1}(\R^d,w_1)^\Theta
\, .
\]
Here we proceed as in Step 3 of the proof of Corollary \ref{interc} taking into account Lemma \ref{interq},
and the formulas (\ref{eq-28b}), (\ref{eq-28c}).
\epr


\section{Complex interpolation of weighted Besov and
Lizorkin-Triebel spaces}
\label{Complex1}


Now we transfer our results on Calder{\'o}n products into results on complex interpolation.


\subsection{Complex interpolation of the spaces $a^s_{p,q}(\Rd,w)$}


We have to take into  account the following supplement to
Proposition \ref{Thm4}.

\begin{Lem}\label{sl1}
Let $X_0,X_1$ be a pair of quasi-Banach sequence lattices.
Then, if both $X_0$ and $X_1$ are analytically convex and at least one is separable, it follows that
$X_0+X_1$ is analytically convex and
\be\label{sl2}
[X_0,X_1]_\Theta = X_0^{1-\Theta}\, X_1^\Theta \, ,\qquad  0<\Theta<1.
\ee
\end{Lem}

\noindent
\bpr
This is the contents of the Remark in front of Theorem 7.10 in \cite{KMM}, see also
\cite{MM}.
\epr

We only need to summarize what we did before.

\begin{Cor}\label{inter}
Let $0 < \Theta < 1$,  $0 < p_0,p_1 \le \infty $, $0 < q_0,q_1 \le \infty$ and $s_0,s_1\in \R$.
Let $p,q$ and $s$ be as in (\ref{eq:inter:1}).
Let $w_0, w_1\in {\mathcal A}_\infty^{\ell oc}$  and define $w$ by the formula
(\ref{eq:inter:2}).
\\
{\rm (i)} If $\max (p_0,p_1)<\infty$ and $\min(q_0,q_1)<\infty$, then
\[
[f^{s_0}_{p_0,q_0}(\R^d,w_0), \, f^{s_1}_{p_1,q_1}(\R^d,w_1)]_\Theta =
f^{s}_{p,q}(\R^d,w)
\]
holds in the sense of equivalent quasi-norms.\\
{\rm (ii)} If $\max(p_0,q_0)<\infty$ or $\max(p_1,q_1)<\infty$, then
\[
[b^{s_0}_{p_0,q_0}(\R^d,w_0), \, b^{s_1}_{p_1,q_1}(\R^d,w_1)]_\Theta =
b^{s}_{p,q}(\R^d,w)
\]
holds in the sense of equivalent quasi-norms.
\end{Cor}

\noindent
\bpr
{\em Step 1.} Proof of (i).
We have to combine Theorem \ref{intera},  Lemma \ref{convex} and Lemma \ref{sl1}.
Because of
$f^s_{p,q} (\Rd,w)$ is separable if, and only if,
$q < \infty$ the claim follows under the condition $\max (q_0,q_1)<\infty$.
Taking into account Lemma \ref{sl1} we replace  $\max (q_0,q_1)<\infty$ by
$\min (q_0,q_1)<\infty$.\\
{\em Step 2.} Proof of (ii).
This time we combine Corollary \ref{interu},   Lemma \ref{convex} and Lemma \ref{sl1}.
Because of
$b^s_{p,q} (\Rd,w)$ is separable if, and only if,
$\max (p,q) < \infty$ the claim follows under the condition $\max (p_0,p_1,q_0,q_1)<\infty$.
By means of Lemma \ref{sl1} we may replace  $\max (p_0,p_1,q_0,q_1)<\infty$ by
$\max(p_0,q_0)<\infty$ or $\max(p_1,q_1)<\infty$.
\epr

\begin{Rem}
 \rm
In case $w_0 = w_1 \equiv 1$ this has been proved in Frazier, Jawerth \cite{FJ2}
($f$-case), Mendez, Mitrea \cite{MM} ($b$-case with $s_0 \neq s_1$)
and Kalton, Mayboroda and Mitrea \cite{KMM} ($b$-case).
The formula
\[
[f^{s_0}_{p_0,q_0}(\R^d,\mu), \, f^{s_1}_{p_1,q_1}(\R^d,\mu)]_\Theta =
f^{s}_{p,q}(\R^d,\mu)
\]
with $\mu$ being a doubling measure has been established by Bownik \cite{Bow2}.
In case $w_0 = w_1 =w \in {\mathcal A}_\infty^{\ell oc}$
\[
[f^{s_0}_{p_0,q_0}(\R^d,w_0), \, f^{s_1}_{p_1,q_1}(\R^d,w_1)]_\Theta =
f^{s}_{p,q}(\R^d,w)
\]
has been proved in Wojciechowska \cite{AW2}.
\end{Rem}

The same type of arguments, this time applied with Theorem \ref{interx} instead of
Theorem \ref{intera} and with
Theorem \ref{intery} instead of
Corollary \ref{interu} yields the next interesting result.
Observe, that all spaces $\mathring{a}^s_{p,q} (\Rd,w)$, $a\in \{b,f\}$, are separable and analytically convex
(use Lemma \ref{convex}).

\begin{Cor}\label{intere}
Let $0 < \Theta < 1$,  $0 < p_0,p_1 \le \infty $, $0 < q_0,q_1 \le \infty$ and $s_0,s_1\in \R$.
Let $p,q$ and $s$ be as in (\ref{eq:inter:1}).
Let $w_0, w_1\in {\mathcal A}_\infty^{\ell oc}$  and define $w$ by the formula
(\ref{eq:inter:2}).
\\
{\rm (i)} If $\max (p_0,p_1)<\infty$, then
\beqq
[\mathring{f}^{s_0}_{p_0,q_0}(\R^d,w_0), \, \mathring{f}^{s_1}_{p_1,q_1}(\R^d,w_1)]_\Theta  & = &
[\mathring{f}^{s_0}_{p_0,q_0}(\R^d,w_0), \, {f}^{s_1}_{p_1,q_1}(\R^d,w_1)]_\Theta
\\
&=&
[{f}^{s_0}_{p_0,q_0}(\R^d,w_0), \, \mathring{f}^{s_1}_{p_1,q_1}(\R^d,w_1)]_\Theta
 = \mathring{f}^{s}_{p,q}(\R^d,w)
\eeqq
holds in the sense of equivalent quasi-norms.\\
{\rm (ii)} Always we have
\beqq
[ \mathring{b}^{s_0}_{p_0,q_0}(\R^d,w_0), \,  \mathring{b}^{s_1}_{p_1,q_1}(\R^d,w_1)]_\Theta & = &
[ \mathring{b}^{s_0}_{p_0,q_0}(\R^d,w_0), \,  {b}^{s_1}_{p_1,q_1}(\R^d,w_1)]_\Theta
\\
&=&
[ {b}^{s_0}_{p_0,q_0}(\R^d,w_0), \,  \mathring{b}^{s_1}_{p_1,q_1}(\R^d,w_1)]_\Theta
=  \mathring{b}^{s}_{p,q}(\R^d,w)
\eeqq
in the sense of equivalent quasi-norms.
\end{Cor}

\begin{Rem}
 \rm
In case $w_0 = w_1 \equiv 1$, Corollary \ref{intere} can be found in
Yang, Yuan and Zhuo \cite{DWC}.
\end{Rem}


\subsection{Complex interpolation of weighted  Besov and Lizorkin-Triebel spaces}
\label{main}


Now we turn to the complex interpolation of the distribution spaces
$F^{s}_{p,q}(\R^d,w)$ and $B^{s}_{p,q}(\R^d,w)$.
For a definition of these classes we refer to
the Appendix.
Observe that neither $F^{s}_{p,q}(\R^d,w)$ nor $B^{s}_{p,q}(\R^d,w)$ are quasi-Banach lattices in general.
Here our main result is as follows.

\begin{T}\label{wichtig}
Let $0 < p_0,p_1 \le \infty $, $0 < q_0,q_1 \le \infty$, $s_0,s_1\in \R$, $0 < \Theta < 1$   and define $p,q$ and
$s$ according to (\ref{eq:inter:1}).
Let $w_0, w_1\in {\mathcal A}_\infty^{\ell oc}$ be local Muckenhoupt weights and $w$ defined as in (\ref{eq:inter:2}).
\\
{\rm (i)} If $\max (p_0,p_1)<\infty$ and $\min(q_0,q_1)<\infty$, then
\[
[F^{s_0}_{p_0,q_0}(\R^d,w_{0}), \, F^{s_1}_{p_1,q_1}(\R^d,w_{1})]_\Theta =
F^{s}_{p,q}(\R^d,w)
\]
holds in the sense of equivalent quasi-norms.\\
{\rm (ii)} If $\max(p_0,q_0)<\infty$ or $\max(p_1,q_1)<\infty$, then
\[
[B^{s_0}_{p_0,q_0}(\R^d,w_{0}), \, B^{s_1}_{p_1,q_1}(\R^d,w_{1})]_\Theta =
B^{s}_{p,q}(\R^d,w)
\]
holds in the sense of equivalent quasi-norms.
\end{T}

\noindent
\bpr
It is enough to combine Corollary \ref{inter} and Proposition \ref{sawano} (in the Appendix).
If either  $p_0 = \infty$ or $p_1 = \infty$, then one has to take into account also Remark \ref{sawano2}.
\epr

\begin{Rem}
 \rm
The above results complement the knowledge on interpolation of weighted  Besov and
Lizorkin-Triebel spaces with Muckenhoupt weights.
Bownik \cite{Bow2} has proved
\[
[F^{s_0}_{p_0,q_0}(\R^d,\mu), \, F^{s_1}_{p_1,q_1}(\R^d,\mu)]_\Theta =
F^{s}_{p,q}(\R^d,\mu)\, ,
\]
where $\mu$ is a doubling measure.
Furthermore, Wojciechowska \cite{AW2} recently  proved
\[
[F^{s_0}_{p_0,q_0}(\R^d, w), \, F^{s_1}_{p_1,q_1}(\R^d,w)]_\Theta =
F^{s}_{p,q}(\R^d,w)\, ,
\]
where $w \in \ca_\infty^{\ell oc}$.
For various interpolation formulas for the real method we refer to Bui \cite{bui-1}
($w \in \ca_\infty$) and Rychkov \cite{Ry} ($w \in \ca_\infty^{\ell oc}$).
\end{Rem}

We finish this subsection by formulating a consequence of Corollary \ref{intere}.
Let $\mathring{A}^{s}_{p,q}(\R^d,w)$ denote the closure of the test functions in
${A}^{s}_{p,q}(\R^d,w)$.
Using the same arguments as above, but replacing Proposition \ref{sawano} by
Proposition \ref{sawano3}, we obtain the following.

\begin{T}\label{interf}
Let $0 < \Theta < 1$,  $0 < p_0,p_1 \le \infty $, $0 < q_0,q_1 \le \infty$ and $s_0,s_1\in \R$.
Let $p,q$ and $s$ be as in (\ref{eq:inter:1}).
Let $w_0, w_1\in {\mathcal A}_\infty^{\ell oc}$  and define $w$ by the formula
(\ref{eq:inter:2}).
\\
{\rm (i)} If $\max (p_0,p_1)<\infty$, then
\beqq
[\mathring{F}^{s_0}_{p_0,q_0}(\R^d,w_0), \, \mathring{F}^{s_1}_{p_1,q_1}(\R^d,w_1)]_\Theta  & = &
[\mathring{F}^{s_0}_{p_0,q_0}(\R^d,w_0), \, {F}^{s_1}_{p_1,q_1}(\R^d,w_1)]_\Theta
\\
&=&
[{F}^{s_0}_{p_0,q_0}(\R^d,w_0), \, \mathring{F}^{s_1}_{p_1,q_1}(\R^d,w_1)]_\Theta
 = \mathring{F}^{s}_{p,q}(\R^d,w)
\eeqq
holds in the sense of equivalent quasi-norms.\\
{\rm (ii)} Always we have
\beqq
[ \mathring{B}^{s_0}_{p_0,q_0}(\R^d,w_0), \,  \mathring{B}^{s_1}_{p_1,q_1}(\R^d,w_1)]_\Theta & = &
[ \mathring{B}^{s_0}_{p_0,q_0}(\R^d,w_0), \,  {B}^{s_1}_{p_1,q_1}(\R^d,w_1)]_\Theta
\\
&=&
[ {B}^{s_0}_{p_0,q_0}(\R^d,w_0), \,  \mathring{B}^{s_1}_{p_1,q_1}(\R^d,w_1)]_\Theta
=  \mathring{B}^{s}_{p,q}(\R^d,w)
\eeqq
in the sense of equivalent quasi-norms.
\end{T}

\begin{Rem}
 \rm
(i) Of course, Theorem \ref{interf} is only of interest in the cases
$\max (p_0,p_1,q_0,q_1)= \infty$.
All other cases are covered by Theorem \ref{wichtig}.\\
(ii) Theorem \ref{interf} has been known since a long time for Besov spaces in the case
$w_0=w_1 \equiv 1 $,  $ 1 <p_0,p_1< \infty$.
We refer to Triebel \cite[Remark~2.4.1/3]{TrI}.
Theorem \ref{interf} in case  $w_0=w_1 \equiv 1 $ and arbitrary $p$'s and $q$'s
has been proved recently in Yang, Yuan, Zhuo \cite{DWC}.

\end{Rem}


\subsection{Shrinking some gaps}
\label{gap}


The results in Subsection \ref{main} do not cover all possible interpolation couples  of the form
\[
[{A}^{s_0}_{p_0,q_0}(\R^d,w_0), \, {A}^{s_1}_{p_1,q_1}(\R^d,w_1)]_\Theta
\, , \qquad A \in \{B,F\}\, .
\]
Those cases, where both spaces are not separable, are not covered.
This means that  a description of
\beqq
[{F}^{s_0}_{p_0,\infty}(\R^d,w_0), \,  {F}^{s_1}_{p_1,\infty}(\R^d,w_1)]_\Theta
& \qquad & 0 < p_0,p_1 < \infty\, ,
\\
 {[}{B}^{s_0}_{p_0,\infty}(\R^d,w_0), \,  {B}^{s_1}_{p_1,\infty}(\R^d,w_1)]_\Theta
& \qquad & 0 < p_0,p_1 \le  \infty\, ,
\\
{[}{B}^{s_0}_{\infty,q_0}(\R^d,w_0), \,  {B}^{s_1}_{\infty ,q_1}(\R^d,w_1)]_\Theta
& \qquad & 0 < q_0,q_1 \le  \infty\, , \quad \min (q_0,q_1)<\infty\, ,
\eeqq
is still open.
We can not fill this gap. However, we can make it smaller.
The method we will apply is based on a result of Shestakov \cite{Sh}
(see also \cite[Remark~4.3.5,~pp.~557]{BK}):
for Banach lattices $X_0, X_1$ and $0 < \Theta < 1$ we have the identity
\be\label{prima}
[X_0,X_1]_ \Theta = \overline{X_0 \cap X_1}^{Y_\Theta}\, , \qquad
Y_\Theta := X_0^{1-\Theta} \, X_1^\Theta\, .
\ee
To use these results we have to switch to sequence spaces for a moment.
Theorem \ref{intera} and Corollary \ref{interu} combined with (\ref{prima}) yield
\[
[{a}^{s_0}_{p_0,q_0}(\R^d,w_0), \, {a}^{s_1}_{p_1,q_1}(\R^d,w_1)]_\Theta
\hookrightarrow {a}^{s}_{p,q}(\R^d,w)\, , \qquad a \in \{b,f\}\, ,
\]
since
\[
\overline{{a}^{s_0}_{p_0,q_0}(\R^d,w_0) \,
\cap \,  {a}^{s_1}_{p_1,q_1}(\R^d,w_1)}^{{a}^{s}_{p,q}(\R^d,w)}
\hookrightarrow {a}^{s}_{p,q}(\R^d,w)\, .
\]
An application of the trivial embedding
\[
[\mathring{A}^{s_0}_{p_0,q_0}(\R^d,w_0), \, \mathring{A}^{s_1}_{p_1,q_1}(\R^d,w_1)]_\Theta
\hookrightarrow
[{A}^{s_0}_{p_0,q_0}(\R^d,w_0), \, {A}^{s_1}_{p_1,q_1}(\R^d,w_1)]_\Theta
\, , \qquad A \in \{B,F\}\, ,
\]
and of Proposition \ref{sawano} yield the following lemma.

\begin{Lem}
\label{einbettung}
Let $1\le p_0,p_1 \le \infty $, $1\le  q_0,q_1 \le \infty$, $s_0,s_1\in \R$, $0 < \Theta < 1$   and define $p,q$ and
$s$ according to (\ref{eq:inter:1}).
Let $w_0, w_1\in {\mathcal A}_\infty^{\ell oc}$ be local Muckenhoupt weights and $w$ defined as in (\ref{eq:inter:2}).
\\
{\rm (i)} If $\max (p_0,p_1)<\infty$, then
\[
\mathring{F}^{s}_{p,q}(\R^d,w)\hookrightarrow
[F^{s_0}_{p_0,q_0}(\R^d,w_{0}), \, F^{s_1}_{p_1,q_1}(\R^d,w_{1})]_\Theta \hookrightarrow
F^{s}_{p,q}(\R^d,w)
\]
holds.\\
{\rm (ii)} Always we have
\[
\mathring{B}^{s}_{p,q}(\R^d,w)\hookrightarrow
[B^{s_0}_{p_0,q_0}(\R^d,w_{0}), \, B^{s_1}_{p_1,q_1}(\R^d,w_{1})]_\Theta \hookrightarrow
B^{s}_{p,q}(\R^d,w)\, .
\]
\end{Lem}

\begin{Rem}
 \rm
Let $w\in \ca_\infty^{\ell oc}$.
The restriction to Banach spaces in Lemma \ref{einbettung} is caused by the use of (\ref{prima}).
Recently, Yang, Yuan and Zhuo \cite{DWC} have proved
\[
[{a}^{s_0}_{p_0,q_0}(\R^d,w), \, {a}^{s_1}_{p_1,q_1}(\R^d,w)]_\Theta
\hookrightarrow {a}^{s}_{p,q}(\R^d,w)\, , \qquad a \in \{b,f\}\, ,
\]
without restrictions on the parameters.
In fact, they did it in the unweighted case, but it can be immediately lifted to the present situation.
Hence, Lemma \ref{einbettung} remains true for values of $p_0,p_1,q_0,q_1 <1$ if restricted to the case of one weight.
\end{Rem}

As mentioned above, in some situations we can go one step further.

\begin{T}
\label{wichtig2}
Let   $s_0,s_1\in \R$ and  $0 < \Theta < 1$.
Further, let either $1 \le  p_0 = p_1 < \infty $ and $s_0 >s_1$
or $1 \le  p_0 < p_1 < \infty $ and
\begin{equation}\label{eq:emb1}
s_0 - \frac{d}{p_0} \ge s_1 - \frac{d}{p_1.}
\end{equation}
As always  $p,q$ and
$s$ are defined according to (\ref{eq:inter:1}).
Let $w \in {\mathcal A}_\infty^{\ell oc}$ be a local Muckenhoupt weight.
Then
\[
\mathring{F}^{s}_{p,\infty}(\R^d,w) = [F^{s_0}_{p_0,\infty}(\R^d,w), \, F^{s_1}_{p_1,\infty}(\R^d,w)]_\Theta \, .
\]
holds.
\end{T}

\noindent
\bpr
The conditions guarantee
\[
F^{s_0}_{p_0,\infty}(\R^d) \hookrightarrow
F^{s}_{p,1}(\R^d)\,
\hookrightarrow
F^{s_1}_{p_1,\infty}(\R^d)\, ,
\]
see \cite[Theorem~2.7.1]{Tr83}.
By means of Proposition \ref{sawano} this yields
\[f^{s_0}_{p_0,\infty}(\R^d)
\hookrightarrow f^{s}_{p,1}(\R^d)\hookrightarrow
f^{s_1}_{p_1,\infty}(\R^d)\, ,
\]
which can be lifted to the weighted case
\[
f^{s_0}_{p_0,\infty}(\R^d,w) \hookrightarrow
f^{s}_{p,1}(\R^d,w) \hookrightarrow
f^{s_1}_{p_1,\infty}(\R^d,w)\,
\]
by using an appropriate isomorphism.
This implies that
\[
f^{s_0}_{p_0,\infty}(\R^d,w) \, \cap f^{s_1}_{p_1,\infty}(\R^d,w) =
f^{s_0}_{p_0,\infty}(\R^d,w)\, .
\]
We claim
\[
\overline{f^{s_0}_{p_0,\infty}(\R^d,w)}^{f^{s}_{p,\infty}(\R^d,w)} =
\mathring{f}^{s}_{p,\infty}(\R^d,w)\, .
\]
To prove this, observe $\mathring{f}^{s}_{p,1}(\R^d,w)\, = f^{s}_{p,1}(\R^d,w)$.
Hence, finite sequences are dense in the set $f^{s_0}_{p_0,\infty}(\R^d,w)$ when equipped with the  quasi-norm of $f^{s}_{p,1}(\R^d,w)$. 
From the trivial embedding $f^{s}_{p,1}(\R^d,w) \hookrightarrow f^{s}_{p,\infty}(\R^d,w)$
we derive the density of the finite sequences in the set
$f^{s_0}_{p_0,\infty}(\R^d,w)$ when equipped with the  quasi-norm of $f^{s}_{p,\infty}(\R^d,w)$. Shestakov's identity (\ref{prima}) yields
\[
\mathring{f}^{s}_{p,\infty}(\R^d,w) = [f^{s_0}_{p_0,\infty}(\R^d,w), \, f^{s_1}_{p_1,\infty}(\R^d,w)]_\Theta \, .
\]
By means of Proposition \ref{sawano} we complete the proof.
\epr

Arguing as before we can prove the following counterpart for Besov spaces.
For the embedding relations we refer to \cite[2.7.1]{Tr83}.

\begin{T}
\label{wichtig3}
Let  $1 \le  q_0,q_1 \le \infty$, $s_0,s_1\in \R$ and  $0 < \Theta < 1$.
Further, let  
\begin{equation}\label{eq:emb2}
s_0 - \frac{d}{p_0} > s_1 - \frac{d}{p_1}
\end{equation}
and either $1 \le p_0 < p_1 \le \infty$ or
$1 \le p_0 = p_1 <\infty$.
As always  $p,q$ and
$s$ are defined according to (\ref{eq:inter:1}).
Let $w \in {\mathcal A}_\infty^{\ell oc}$ be a local Muckenhoupt weight.
Then
\[
\mathring{B}^{s}_{p,q}(\R^d,w) = [B^{s_0}_{p_0,q_0}(\R^d,w), \, B^{s_1}_{p_1,q_1}(\R^d,w)]_\Theta \, .
\]
holds.
\end{T}

\begin{Rem}
\rm (i) The case $w \equiv 1$, $1 <p_0 = p_1 <\infty$, $q_0 = q_1 =q= \infty$ has been known before,
we refer to Triebel \cite[Theorem~2.4.1]{TrI}.\\
(ii) Let us comment on the conditions \eqref{eq:emb1} and \eqref{eq:emb2}. Let $w \equiv 1, 1\le p_0<p_1 \le\infty$ and
\begin{equation}\label{eq:emb3}
s_0 - \frac{d}{p_0} \le s_1 - \frac{d}{p_1}.
\end{equation}
Furthermore, for $j\in\N_0$ let  $K_j$ be a subset of $\Z^d$ with cardinality
\[
\#K_j=\lceil 2^{-j\{(s_1-s_0)\cdot\frac{1}{1/p_1-1/p_0}-d\}}\rceil\, ,
\]
where $\lceil t\rceil$
denotes the smallest integer larger than or equal to $t\in\R.$ We define a sequence $\lambda :=\{\lambda_{j,k}\}_{j,k}$ by
\[
\lambda_{j,k}=\begin{cases}\displaystyle 2^{j\cdot\frac{p_1s_1-p_0s_0}{p_0-p_1}}& \qquad \text{if}\quad  k\in K_j,\\
0
& \qquad \text{otherwise}.\end{cases}
\]
A simple calculation shows
that $\lambda\in b^s_{p,\infty}(\R^d)\setminus \mathring{b}^s_{p,\infty}(\R^d)$ as well as
$\lambda\in b^{s_0}_{p_0,\infty}(\R^d)\cap b^{s_1}_{p_1,\infty}(\R^d).$ The result of Shestakov \eqref{prima} then yields
\[
\lambda\in \overline{b^{s_0}_{p_0,\infty}(\R^d)\cap b^{s_1}_{p_1,\infty}(\R^d)}^{b^s_{p,\infty}(\R^d)}=[b^{s_0}_{p_0,\infty}(\R^d),b^{s_1}_{p_1,\infty}(\R^d)]_\theta.
\]
Hence,  the embedding $\mathring{b}^{s}_{p,q}(\R^d) \hookrightarrow [b^{s_0}_{p_0,q_0}(\R^d), \, b^{s_1}_{p_1,q_1}(\R^d)]_\Theta$ is strict in
this case. With other words, (\ref{eq:emb2}) is also necessary for the validity of the interpolation formula in
Theorem \ref{wichtig3}.
\end{Rem}


\section{Complex interpolation of radial subspaces of Besov and
Lizorkin-Triebel spaces}
\label{complex}


In a series of papers the authors have studied radial subspaces of Besov and
Lizorkin-Triebel spaces, see \cite{SS, SSV1, SS2}.
The motivation came from the interesting interplay of decay and smoothness properties of radial functions as expressed in its simplest form
in the radial Lemma of Strauss \cite{Strauss} and with important applications for
the compactness of some embeddings.
We refer also to Lions \cite{Lions} and Cho and Ozawa \cite{CO} in this connection. Let us recall the convention, that if $E$ denotes a space of functions on $\Rd$ then by $RE$ we mean the subset of radial functions
in $E$ and we endow this subset with the same quasi-norm as the original space.


\subsection{The main result for radial subspaces}


In \cite{LSR} one of the authors has proved that in case $p,q\ge 1$
the spaces  $RB^{s}_{p,q} (\Rd), \, RF^{s}_{p,q} (\Rd) $ are  complemented
subspaces of $\bspq (\Rd)$ and $\fspq (\Rd)$, respectively.
By means of the method of retraction and coretraction, see, e.g.,
Theorem 1.2.4 in \cite{TrI}, this allows to transfer the interpolation formulas
for the original spaces $B^{s}_{p,q} (\Rd), \, F^{s}_{p,q} (\Rd) $
to their radial subspaces.
\\
Such simple arguments do not seem to be available in case of quasi-Banach spaces.
In \cite{SSV1} we announced the following result concerning complex interpolation
of radial subspaces of Besov and Lizorkin-Triebel spaces.

\begin{T}\label{interpol3}
Let  $0 <  p_0 ,  p_1 < \infty$, $0 <  q_0 ,\,  q_1 \le \infty$, $s_0, s_1 \in \R$, and $0 < \Theta < 1$.
Define $s:= (1-\Theta)\, s_0 + \Theta \, s_1$,
\[
\frac 1p := \frac{1-\Theta}{p_0} + \frac{\Theta}{p_1}
\qquad \mbox{and}\qquad \frac 1q := \frac{1-\Theta}{q_0} +
\frac{\Theta}{q_1} \, .
\]
{\rm (i)} Let $\min (q_0,q_1) <\infty$.
Then we have
\[
RB^s_{p,q} (\Rd)  =   \Big[RB^{s_0}_{p_0,q_0} (\Rd),
RB^{s_1}_{p_1,q_1} (\Rd)\Big]_{\Theta} \, .
\]
{\rm (ii)} Let $\min (q_0,q_1) <\infty$.
Then we have
\[
RF^s_{p,q} (\Rd) =   \Big[RF^{s_0}_{p_0,q_0} (\Rd),
RF^{s_1}_{p_1,q_1} (\Rd)\Big]_{\Theta} \, .
\]
\end{T}



\subsection{The proof of Theorem 29}


We are going to prove Theorem \ref{interpol3}.
Our main idea consists in a reduction of this problem to the weighted case.
Therefore we need to recall some results from \cite{SSV1}.
\\
For a real number $s$ we denote by $[s]$ the integer part, i.e. the largest
integer $m$ such that $m \le s$. We put
$w_{d-1}(t):= |t|^{d-1}$, $t\in \R$, $d\ge 2$. Finally, if $f:\Rd\to\C$ is a radial function, we denote by
$$
(\tr f)(t)=f(t,0,\dots,0),\quad t\in \R
$$
its trace. The corresponding extension operator is then given by
$$
(\ext g)(x)=g(|x|),\quad x\in\R^d,
$$
where $g$ is an even function on $\R$.

\begin{Prop}\label{main5}
Let $d \ge 2$, $0 < p< \infty$ and $0 < q \le \infty$.
\\
{\rm (i)} Suppose either  $s > d(\frac 1p - \frac 1d)$
or $s= d(\frac 1p - \frac 1d) $ and $q\le 1$.
Then the  mapping $\tr $
is a linear isomorphism of $R\bspq (\Rd)$ onto  $R\bspq (\R, w_{d-1})$
with inverse $\ext$.
\\
{\rm (ii)} Suppose  either  $s > d(\frac 1p - \frac 1d)$
or $s = d(\frac 1p - \frac 1d)$ and $0< p \le 1$.
Then the mapping $\tr $
is a linear  isomorphism of  $R\fspq (\Rd)$ onto  $R\fspq (\R, w_{d-1})$
with inverse  $\ext$.
\end{Prop}

\begin{Rem}
 \rm
Let $p>1$.
As it is well known, the weight $w(x):= |x|^\alpha$ belongs to $\ca_p (\R)$
if, and only if, $-1 < \alpha <  p-1$.
This implies
$w_{d-1} \in {\mathcal A}_p (\R)$ for any $p>d$, see \cite[pp.~218]{stein}.
\end{Rem}

We would like to apply Theorem \ref{wichtig} with respect to $R\bspq (\R, w_{d-1})$ and $R\fspq (\R, w_{d-1})$, respectively.
Therefore we consider the following mapping:
\[
T f (x) := \frac 12 \, (f(x) + f(-x))\, , \qquad x \in \R\, .
\]
Of course,
\[
T \in \cl (\aspq (\R, w_{d-1}), R\aspq (\R, w_{d-1}))\, , \qquad A \in \{B,F\}\, .
\]
Furthermore, $T f = f$ if $f \in R\aspq (\R, w_{d-1})$, i.e.,
$R\aspq (\R, w_{d-1})$ is a retract of $\aspq (\R, w_{d-1})$.
By the standard method of retraction and coretraction, see \cite[1.2.4]{TrI} and Lemma 7.11 in \cite{KMM}, we obtain the following.

\begin{Lem}\label{gut}
Let $d \ge 2$,  $0 < \Theta < 1$, $0 < p_0,p_1 \le  \infty $, $0 < q_0,q_1 \le \infty$,  $s_0,s_1\in \R$ and define $p,q$ and
$s$ according to (\ref{eq:inter:1}).
\\
{\rm (i)} If $\max (p_0,p_1)<\infty$ and $\min(q_0,q_1)<\infty$, then
\[
[RF^{s_0}_{p_0,q_0}(\R,w_{d-1}), \, RF^{s_1}_{p_1,q_1}(\R,w_{d-1})]_\Theta =
RF^{s}_{p,q}(\R,w_{d-1})
\]
holds in the sense of equivalent quasi-norms.\\
{\rm (ii)} If $\max(p_0,q_0)<\infty$ or $\max(p_1,q_1)<\infty$, then
\[
[RB^{s_0}_{p_0,q_0}(\R,w_{d-1}), \, RB^{s_1}_{p_1,q_1}(\R,w_{d-1})]_\Theta =
RB^{s}_{p,q}(\R,w_{d-1})
\]
holds in the sense of equivalent quasi-norms.
\end{Lem}

The next step consists in combining Lemma \ref{gut} and Proposition \ref{main5}.

\begin{Lem}\label{gut2}
Let $d \ge 2$,  $0 < \Theta < 1$, $0 < p_0,p_1 < \infty $, $0 < q_0,q_1 \le \infty$,  $s_0,s_1\in \R$ and define $p,q$ and
$s$ according to (\ref{eq:inter:1}).
Furthermore, we suppose
\be\label{eq-30b}
s_i > d\, (\frac {1}{p_i} - \frac 1d) \, , \qquad i=0,1\, ,
\ee
If $\min(q_0,q_1)<\infty$, then
\[
[RF^{s_0}_{p_0,q_0}(\R^d), \, RF^{s_1}_{p_1,q_1}(\R^d)]_\Theta =
RF^{s}_{p,q}(\R^d)
\]
and
\[
[RB^{s_0}_{p_0,q_0}(\R^d), \, RB^{s_1}_{p_1,q_1}(\R^d)]_\Theta =
RB^{s}_{p,q}(\R^d)
\]
hold in the sense of equivalent quasi-norms.
\end{Lem}

\noindent
{\bf Proof of Theorem \ref{interpol3}}.
It is enough to remove the restrictions for $s_0$ and $s_1$ in Lemma \ref{gut2}, see (\ref{eq-30b}).
But this is an easy task. Let $\sigma \in \R$.
We consider the family of  lifting operators
\[
I_\sigma f := \cfi [(1+ |\xi|^2)^{\sigma/2} \, \cf f (\xi)] (\, \cdot \, )\, , \qquad f \in \cs' (\Rd)\, .
\]
As it is well known, see, e.g., \cite[Theorem~2.3.8]{Tr83},
$I_\sigma$ is an isomorphism, which maps
$\aspq (\Rd)$ onto $A^{s-\sigma}_{p,q} (\Rd)$.
By standard properties of the Fourier transform we deduce that $f \in R\aspq (\Rd)$ implies
$I_\sigma f \in RA^{s-\sigma}_{p,q} (\Rd)$.
By the same argument, $I_\sigma$ is an isomorphism, which maps
$R\aspq (\Rd)$ onto $RA^{s-\sigma}_{p,q} (\Rd)$.
Hence,
\[
[I_\sigma (RF^{s_0}_{p_0,q_0}(\R^d)), \, I_\sigma (RF^{s_1}_{p_1,q_1}(\R^d))]_\Theta =
I_\sigma (RF^{s}_{p,q}(\R^d)) \, .
\]
Now, Theorem \ref{interpol3} follows from Lemma \ref{gut2} by choosing $\sigma $ appropriate.


\section{Appendix -- Muckenhoupt weights and function spaces}


For convenience of the reader we collect some definitions and properties
around Muckenhoupt weights and associated weighted function spaces.


\subsection{Muckenhoupt weights}


A weight function (or simply a weight) is a nonnegative and measurable  function on $\Rd$.
We collect a few facts including the definition of Muckenhoupt
and local Muckenhoupt weights.
As usual, $p'$ is related to $p$ via the formula  $1/p + 1/p'=1$.

\begin{Def}
Let $1 < p < \infty$.
Let $w$ be a nonnegative, locally integrable function on $\Rd$.
\\
(i) Then $w$ belongs to the  Muckenhoupt class $\ca_p$, if
\[
A_p (w) := \sup_{B} \, \Big(\frac{1}{|B|} \int_B w(x)\, dx\Big)^{1/p} \, \cdot \,
\Big(\frac{1}{|B|} \int_B w(x)^{-p'/p} \, dx\Big)^{1/p'} <\infty\, ,
\]
where the supremum is taken with respect to all balls $B$ in $\Rd$.
\\
(ii) The weight $w$ belongs to the local  Muckenhoupt class $\ca_p^{\ell oc} $,
if we restrict the set of admissible balls in the supremum in (i) to those with volume $\le 1$.
We put
\[
A_p^{\ell oc} (w) := \sup_{|B|\le 1} \, \Big(\frac{1}{|B|} \int_B w(x)\, dx\Big)^{1/p} \, \cdot \,
\Big(\frac{1}{|B|} \int_B w(x)^{-p'/p} \, dx\Big)^{1/p'} <\infty\, .
\]
\end{Def}

The  classes  $\ca_\infty$  and  $\ca_\infty^{\ell oc}$ are  defined as
\[
\ca_\infty := \bigcup_{p>1} \ca_p
\qquad \mbox{and}\qquad \ca_\infty^{\ell oc}  := \bigcup_{p>1} \ca_p^{\ell oc}\, ,
\]
respectively.

\begin{Rem}
 \rm
(i) Good sources for Muckenhoupt weights are Stein \cite{stein}, Garc\'\i a-Cuerva and Rubio de Francia \cite{GC} and the graduate text \cite{Du} of Duoandikoetxea.
\\
(ii) The classes of local Muckenhoupt weights $\ca_p^{\ell oc}$ have been introduced by
Rychkov \cite{Ry}.
\end{Rem}

The following lemma of Rychkov \cite{Ry} will be of some use.

\begin{Lem}\label{local}
 Let $1 <p\le \infty$, $w \in \ca_p^{\ell oc} $, and $Q$ be a cube with sides parallel to the axes and volume $1$.
Then there exists a weight $\overline{w} \in \ca_p $ s.t.
\[
\overline{w} = {w} \qquad \mbox{on $Q$ and}\qquad A_p (\overline{w}) \le c\, A_p^{\ell oc} (w)\, .
\]
Here $c$ does not depend on $w$ and $Q$.
\end{Lem}

By $M f$ we denote the Hardy-Littlewood maximal function of $f$, given by
\[
Mf (x):= \sup_B \frac{1}{|B|} \int_B |f(y)|\, dy\, ,
\]
where the supremum is taken with respect to all balls in $\Rd$ with center $x$.
Furthermore, by $M^{\ell oc}$ we denote the following local counterpart of $M$, namely
\[
M^{\ell oc}f (x):= \sup_Q \frac{1}{|Q|} \int_Q |f(y)|\, dy\, ,
\]
where the supremum is taken with respect to all cubes $Q$ containing $x$, with sides parallel to the axes and volume $\le 1$.
The weighted Lebesgue space $L_p (\Rd,w)$
is the collection of all measurable functions $f: \, \Rd \to \C$ such that
\[
\| \, f \, |L_p (\Rd,w)\| := \Big(\int_{\Rd} |f(x)|^p \, w(x)\, dx \Big)^{1/p}<\infty \, .
\]
In case $p=\infty$ we are back in the unweighted situation, i.e., $w \equiv 1$.
We shall also need the following maximal inequality.
Let $1 < p< \infty$, $1 < q \le \infty$ and $w\in \ca_p$. Then there exists
a constant $c$ such that
\be\label{vachtang}
\Big(\int_{\Rd} \Big(\sum_{j=0}^\infty |M f_j(x)|^q\Big)^{p/q} \, w(x)\, dx \Big)^{1/p}\le c\,
\Big(\int_{\Rd} \Big(\sum_{j=0}^\infty |f_j(x)|^q\Big)^{p/q} \, w(x)\, dx \Big)^{1/p}
\, ,
\ee
holds for all sequences $(f_j)_j \subset L_p (\Rd,w)$, see \cite{Ko}, \cite{AJ} or \cite[Theorem~V.3.1]{stein}.
Rychkov \cite{Ry} proved the local version: for $1 < p< \infty$, $1 < q \le \infty$
and $w\in \ca_p^{\ell oc}$ there exists a constant $c$ s.t.
\be\label{slava}
\Big(\int_{\Rd} \Big(\sum_{j=0}^\infty |M^{\ell oc} f_j(x)|^q\Big)^{p/q} \, w(x)\, dx \Big)^{1/p}\le c\,
\Big(\int_{\Rd} \Big(\sum_{j=0}^\infty |f_j(x)|^q\Big)^{p/q} \, w(x)\, dx \Big)^{1/p}
\, ,
\ee
holds for all sequences $(f_j)_j \subset L_p (\Rd,w)$.
\\
We need one further property of Muckenhoupt weights.

\begin{Lem}\label{stein2}
Let $0 < \Theta < 1$ and  $0 < p_0,p_1 \le \infty $.  We put
\[
\frac 1 p := \frac{1-\Theta}{p_0} + \frac{\Theta}{p_1}\, .
\]
(i) Let $w_0, w_1\in {\ca}_\infty$ and define
\[
w := w_{0}^{\frac{(1-\Theta)p}{p_0}}w_1^{\frac{\Theta p}{p_1}}.
\]
Then $w\in {\ca}_\infty$.
\\
(ii) If $w_0, w_1$ belong to ${\ca}_\infty^{\ell oc}$, then $w\in {\ca}_\infty^{\ell oc}$.
\end{Lem}

\noindent
\bpr
{\em Step 1.} We prove (i).
If $w_0, w_1\in {\ca}_\infty$ then there exist $r_0,r_1 \in (1,\infty)$ such that $w_i \in {\ca}_{r_i}$, $i=0,1$. First, let $\max (p_0,p_1)<\infty$. If $r_i\le p_i$, then by the monotonicity of Muckenhoupt classes $w_i\in {\ca}_{p_i}$  and this implies that $w \in {\ca}_{p}$, see Stein  \cite[V.6.1(a), pp.~218]{stein}.
If $\max\big(\frac{r_0}{p_0},\frac{r_1}{p_1}\big)> 1$, then we choose $\alpha >\max\big(\frac{r_0}{p_0},\frac{r_1}{p_1}\big)$.
Now $w_i\in {\ca}_{\alpha p_i}$,  $i=0,1$, follows because of
\[
\frac{1}{\alpha \, p} = \frac{1-\Theta}{\alpha \, p_0} + \frac{\Theta}{\alpha \, p_1}
\qquad \mbox{and} \qquad w= w_{0}^{\frac{(1-\Theta) \alpha \, p}{\alpha \, p_0}}\, w_1^{\frac{\Theta \, \alpha\,  p}{\alpha \, p_1}} \, .
\]
So  applying the same argument as before we get $w\in {\ca}_{\alpha p} \subset {\ca}_{\infty}$.
For the remaining cases  $\max (p_0,p_1) = \infty$ it is enough to observe that the function
$w(x)= 1$, $x\in \Rd$, belongs to $\ca_\infty$ ($p_0 = p_1= \infty$)
and in case $0 < p_0 < \infty = p_1$ we have $(1-\theta)p/p_0 = 1$, i.e., $w=w_0$.
\\
{\em Step 2}. The monotonicity of the local  Muckenhoupt classes $\ca_p^{\ell oc}$
has been proved in \cite{Ry}. The above used result from
\cite[V.6.1(a), pp.~218]{stein} is based on H\"older's inequality.
For that reason it carries over to the local situation.
Alternatively one could argue with Lemma \ref{local}.
\epr

%


\subsection{Weighted Besov and Lizorkin-Triebel spaces}


Now we introduce weighted  Besov and Lizorkin-Triebel spaces. Since we work with local Muckenhoupt weights  we need larger  space of distributions  than the spaces of tempered distributions. Recall that the class ${\cal A}_\infty^{loc}$ contains weights of exponential growth. We follows the ideas of Rychkov that are based on local reproducing formula, cf. \cite{Ry}. 

Let $S_e(\Rd)$ denote the set of all  $\psi \in C^\infty (\Rd)$ such that  
\[
q_N(\psi) = \sup_{x\in \Rd} e^{N|x|} \sum_{|\alpha|\le N} |\partial^\alpha \psi(x)| \,<\, \infty\qquad \text{for all} \quad N\in \N_0\, .  
\]
Then $S_e(\Rd)$, equipped with the topology generated by the system of semi-norms $q_N$, is a locally convex space. Its dual space $S'_e(\Rd)$ can be identified with a subspace of the space of distributions ${\cal D}'(\Rd)$ in the obvious way. 

Let $\varphi_0 \in C_0^\infty (\Rd)$  and  $\varphi(x) = \varphi_0(x) - 2^{-d} \varphi_0(\frac{x}{2})$ be   functions such that
\[
\int_{\Rd} \varphi_0(x) dx\, \not=\, 0 \qquad \text{and}\qquad \int_{\Rd} x^\beta\varphi(x) dx\, =\, 0 
\]
for any multiindex $\beta\in \N^d_0$, $|\beta|\le B$, where $B$ is a fixed natural number. 
   
\begin{Def}
Let $0 < q \le \infty$, $s\in \R$ and $w \in \ca_\infty^{\ell oc}$. Moreover let  $B>[s]$. \\
{\rm (i)}
Let $0 < p< \infty$. Then the weighted Besov space
$B^s_{p,q}(\Rd,w)$ is the collection of all $f \in S'_e(\Rd)$
such that
\[
\| \, f \, |B^s_{p,q}(\Rd,w)\| := \Big(
\sum_{j=0}^\infty 2^{jsq}\, \| \, \varphi_j * f \, |L_p(\Rd,w)\|^q\Big)^{1/q}<\infty\, .
\]
{\rm (ii)}
Let $0 < p< \infty$. Then the weighted Triebel-Lizorkin space
$F^s_{p,q}(\Rd,w)$ is the collection of all $f \in S'_e(\Rd)$
such that
\[
\| \, f \, |F^s_{p,q}(\Rd,w)\| := \Big\|\, \Big(
\sum_{j=0}^\infty 2^{jsq}\, | \, \varphi_j * f \, |^q \Big)^{1/q} \Big| L_p (\Rd,w)\Big\| <\infty\, .
\]
\end{Def}

\begin{Rem}
\label{unendlich}
\rm
(i)
For $w\equiv 1$ we are in the unweighted case.
The associated spaces are denoted by $B^s_{p,q}(\Rd)$ and $F^s_{p,q}(\Rd)$. The above definition coincides with characterization of $B^s_{p,q}(\Rd)$ and $F^s_{p,q}(\Rd)$ by so called local means,  cf. \cite{Tr92} and \cite{BPT-2}. 
\\
(ii)
Observe, that we did not define weighted spaces with $p=\infty$. However, it will be convenient for us to use
the convention $B^s_{\infty,q}(\Rd,w):= B^s_{\infty,q}(\Rd)$. 
\\
(iii) Weighted Besov and Lizorkin-Triebel spaces with $w\in {\ca}_\infty$ have been first studied systematically by Bui \cite{bui-1,bui-2}, cf. also \cite{BPT-1}  and \cite{BPT-2}.
In addition we refer  to Haroske, Piotrowska \cite{HP} and \cite{HS}.
Standard references for unweighted spaces are
the monograph's \cite{Pe,Tr83,Tr92,Tr06} as well as \cite{FJ2}.
These classes with ${\ca}_\infty^{\ell oc}$ weights
have been treated by Rychkov \cite{Ry},  Izuki and Sawano \cite{IS}, and Wojciechowska \cite{AW,AW2}.
Different types of measure have been considered by
Bownik and Ho \cite{Bow1} and  Bownik \cite{Bow2}.
\end{Rem}

\subsection*{Wavelet characterizations of weighted spaces}
\label{wave}


Here we need the following result.

\begin{Prop}\label{sawano}
Let $0 < p < \infty $,  $0 < q \le \infty$ and $s\in \R$.
Let $w \in \ca_\infty^{\ell oc}$.
Then there exists a linear isomorphism $T$ which maps
$B^s_{p,q} (\Rd,w)$ onto $b^{s+d/2}_{p,q} (\Rd,w)$
and $F^s_{p,q} (\Rd,w)$ onto $f^{s+d/2}_{p,q} (\Rd,w)$.
\end{Prop}

\begin{Rem}\label{sawano2}
 \rm
(i) The mapping $T$ is generated by an appropriate  wavelet system.
A proof of Proposition \ref{sawano} can be found in Izuki and Sawano \cite{IS},
cf. also  Wojciechowska \cite{AW} for a different proof.
In case of Besov spaces and $w \in \ca_\infty$ we also  refer  to \cite{HS} and Bownik, Ho \cite{Bow1} in this context.
\\
(ii) Proposition \ref{sawano} extends to $p=\infty$ for Besov spaces, see Remark \ref{unendlich}.
Wavelet characterizations of unweighted Besov spaces are proved at various places,
we refer to Meyer \cite{Me}, Kahane and Lemarie-Rieusset \cite{KL}, Triebel \cite[3.1.3]{Tr06} and Wojtaszczyk \cite{woj}.
\end{Rem}

There is a little supplement to the previous proposition
dealing with the classes $\mathring{A}^s_{p,q} (\Rd,w)$ and $\mathring{a}^s_{p,q} (\Rd,w)$
($A \in \{B,F\}$, $a \in \{b,f\}$).

\begin{Prop}\label{sawano3}
Let $0 < q \le \infty$,  $s\in \R$ and $w \in \ca_\infty^{\ell oc}$.
\\
(i) Let $0 < p < \infty $.
Then there exists a linear isomorphism $T$ which maps
$\mathring{B}^s_{p,q} (\Rd,w)$ onto $\mathring{b}^{s+d/2}_{p,q} (\Rd,w)$
and $\mathring{F}^s_{p,q} (\Rd,w)$ onto $\mathring{f}^{s+d/2}_{p,q} (\Rd,w)$.
\\
(ii) There exists a linear isomorphism $T$ which maps
$\mathring{B}^s_{\infty,q} (\Rd,w)$ onto $\mathring{b}^{s+d/2}_{\infty,q} (\Rd,w)$.
\end{Prop}

The proof of Proposition \ref{sawano3} follows the same pattern as the proof of Proposition \ref{sawano}. We leave out the details but see \cite[Proposition~2.1]{DWC} for the unweighted case.


\end{document}